\documentclass[11pt,a4paper]{amsproc}
\usepackage[T1]{fontenc}
\usepackage{lmodern}
\usepackage[a4paper,left=22mm,right=22mm,top=30mm,bottom=22mm,headheight=14pt,headsep=6mm,footskip=10mm]{geometry}
\usepackage{amsmath,amssymb,mathtools}
\usepackage{microtype}
\usepackage{needspace}
\usepackage{etoolbox}
\makeatletter
\patchcmd{\@maketitle}{\global\topskip8pc\relax}{\global\topskip1.5pc\relax}{}{\PackageError{a4layout}{Title spacing patch failed}{}}
\makeatother
\usepackage[hidelinks]{hyperref}
\hypersetup{pdftitle={The metric theory of tracial matrix algebras does not converge},pdfauthor={Vadim Alekseev, Andreas Thom}}

\newtheorem{maintheorem}{Theorem}

\newtheorem{isomorphismtheorem}{Theorem}

\newtheorem{nonisomorphismtheorem}{Theorem}

\newtheorem{theorem}{Theorem}[section]
\newtheorem{lemma}[theorem]{Lemma}
\newtheorem{proposition}[theorem]{Proposition}
\newtheorem{corollary}[theorem]{Corollary}
\theoremstyle{definition}
\newtheorem{definition}[theorem]{Definition}
\theoremstyle{remark}

\numberwithin{equation}{section}

\newcommand{\C}{\mathbb C}
\newcommand{\N}{\mathbb N}
\newcommand{\M}{\mathcal M}
\newcommand{\U}{\mathcal U}
\newcommand{\V}{\mathcal V}
\newcommand{\dist}{\operatorname{dist}}
\newcommand{\Sym}{\operatorname{Sym}}
\newcommand{\Th}{\operatorname{Th}}
\newcommand{\dminus}{\mathbin{\dot{-}}}
\allowdisplaybreaks[1]
\title[The metric theory of tracial matrix algebras does not converge]{The metric theory of tracial matrix algebras\\does not converge}
\date{5 November 2026}
\author[V. Alekseev]{Vadim Alekseev}
\author[A. Thom]{Andreas Thom}
\address{Institut f\"ur Geometrie, Technische Universit\"at Dresden, 01062 Dresden, Germany}
\email{vadim.alekseev@tu-dresden.de, andreas.thom@tu-dresden.de}

\begin{document}
\begin{abstract}
We show that the theories of tracial matrix algebras and metric
unitary groups do not converge with the dimension. In the main argument, we recover the
canonical diagonal using spectral gap and internality. Uniform
definability of the diagonal and its normalizer then transfers the
distinct theories of universal sofic groups to matrix ultraproducts. The same methods show that isomorphisms preserve every
internal subalgebra and give fixed finite distance formulas over
six unitaries. Finally, we
construct an explicit continuum family of matrix ultraproducts with
distinct theories and trivial first-order fundamental groups.
\end{abstract}
\maketitle
\tableofcontents

\section{Introduction}\label{sec:introduction}

We study the complete continuous first-order theories of the matrix
algebras \(M_n(\C)\), equipped with the normalized trace
\(\tau_n=n^{-1}\operatorname{Tr}\) and the metric induced by
\(\|x\|_2=\tau_n(x^*x)^{1/2}\) on operator-norm bounded balls.
A sentence has a real value in each algebra, and two tracial algebras
are elementarily equivalent if all these values agree.
Farah, Hart and Sherman asked whether every sentence value converges
as \(n\to\infty\) \cite[Question~6.6]{FHS2}. To express this
question in terms of ultraproducts, we use \L{}o\'s' theorem:
\[
 \varphi^{\prod_{\mathcal U}M_n(\C)}
 =\lim_{\mathcal U}\varphi^{M_n(\C)},
\]
so convergence is equivalent to elementary equivalence of all
nonprincipal tracial ultraproducts of this sequence. Our first
result gives a negative answer.

\begin{maintheorem}\label{cor:many-theories}
There are \(2^{\aleph_0}\) pairwise non-elementarily equivalent
tracial matrix ultraproducts of the sequence \((M_n(\mathbb C),\tau_n)\),
and the same holds for metric ultraproducts of \((U(n),d_2)\) in
the pure metric group language, where \(d_2(u,v)=\|u-v\|_2\).
Consequently, the sequences of complete continuous theories
\[
 \Th(M_n(\mathbb C),\tau_n),\qquad \Th(U(n),d_2)
\]
do not converge as \(n\to\infty\).
\end{maintheorem}

Purely universal and purely existential sentences do converge,
since every tracial matrix ultraproduct has the universal and
existential theories of the hyperfinite \(\mathrm{II}_1\) factor
\cite[Remark~6.5(1), equation~(6.2)]{FHS3}.

To distinguish the theories, we recover the diagonal and use its
unitary normalizer. The quotient by the diagonal unitaries is a
universal sofic group, so uniform recovery of this quotient allows
us to apply the nonconvergence result of \cite{AT} to matrix algebras.
We first fix the class of ultraproducts and the notion of internality
used in this recovery.

\begin{definition}\label{def:tracial-matrix-ultraproduct}
Let \(\mathcal U\) be a nonprincipal ultrafilter on \(\N\), and
let \(\mathbf d=(d_n)\) be a sequence of positive integers with
\(d_n\to\infty\). We call
\begin{equation}\label{eq:Mud}
 \M_{\mathcal U,\mathbf d}
 =\prod_{\mathcal U}(M_{d_n}(\C),\tau_{d_n})
\end{equation}
a \emph{tracial matrix ultraproduct}. Its elements are classes
\([x_n]_{\mathcal U}\) of operator-norm bounded sequences, modulo
sequences whose \(2\)-norms tend to zero along \(\mathcal U\).
The operations are induced coordinatewise, and
\(\tau([x_n]_{\mathcal U})=\lim_{\mathcal U}\tau_{d_n}(x_n)\).
A unital von Neumann subalgebra \(A\subseteq\M_{\mathcal U,\mathbf d}\)
is \emph{internal} if \(A=[A_n]_{\mathcal U}\) for unital
\(*\)-subalgebras \(A_n\subseteq M_{d_n}(\C)\).
An \emph{internal MASA} is an internal abelian subalgebra that is
maximal among internal abelian subalgebras.
\end{definition}

Throughout, \(\M_{\mathcal U,\mathbf d}\) and
\(\M_{\mathcal V,\mathbf e}\) denote tracial matrix ultraproducts
in this sense. To recognize the diagonal, we will use
Lemma~\ref{lem:commutant}: internal MASAs are maximal abelian
in the ambient algebra and are unitarily
conjugate to the canonical diagonal \([D_{d_n}]_{\mathcal U}\),
where \(D_d\subseteq M_d(\C)\) denotes the diagonal algebra.
Isomorphisms are complex \(*\)-isomorphisms; they preserve the trace
because these ultraproducts are \(\mathrm{II}_1\) factors.

To use this recovery for elementary equivalence, we work in the
continuous language of tracial \(*\)-algebras on operator-norm
bounded sorts. A predicate is \emph{definable} if it is a uniform
limit of formulas, and a nonempty closed set is definable if its
distance predicate is definable. We allow parameters when specified,
use bounded real ranges, and write
\(s\dminus t=\max\{s-t,0\}\). See
\cite[Section~3.2]{FHLRTVW} for quantification over definable sets.
The following theorem gives a uniform way to recover an internal
diagonal from finitely many parameters.

\begin{maintheorem}\label{thm:diagonal-code}
There are fixed finite formulas \(\Delta(\mathbf z)\) and
\(\delta(x;\mathbf z)\), with \(\mathbf z\) a tuple of six
contractions, such that in every tracial matrix ultraproduct
\(\M_{\mathcal U,\mathbf d}\):
\begin{enumerate}
\item the zero-set of \(\Delta\) is nonempty;
\item every \(\mathbf z\) in this zero-set defines an internal
MASA \(D_{\mathbf z}\) with
\[
 \delta(x;\mathbf z)=\dist_2(x,(D_{\mathbf z})_{\leq1})
 \qquad(x\in(\M_{\mathcal U,\mathbf d})_{\leq1}).
\]
\end{enumerate}
Every such MASA is unitarily conjugate to the canonical diagonal.
\end{maintheorem}

The diagonal cannot be the commutant of a finite tuple
\cite[Theorem~2.1(1)]{PopaIndependence}, so we recover it as the
relative join of two abelian commutants with spectral gap. The
spectral gap internality theorem of Liu and the authors
\cite[Theorem~B]{ALT} makes these commutants internal. Uniform
abelian Hilbert--Schmidt stability then ensures that their relative
join is internal as well.

The finite code and conjugacy of internal MASAs will allow us to
recover the theory of the pair \((\M_{\mathcal U,\mathbf d},D)\)
from that of the ambient algebra. Here the pair is the expansion by
the distance to the operator-norm unit ball of \(D\). We obtain
the following reduction.

\begin{maintheorem}\label{thm:reduction}
Let \(D=[D_{d_n}]_{\mathcal U}\subseteq\M_{\mathcal U,\mathbf d}\)
and \(F=[D_{e_n}]_{\mathcal V}\subseteq\M_{\mathcal V,\mathbf e}\)
be the canonical diagonal subalgebras of two tracial matrix
ultraproducts. The following hold.
\begin{enumerate}
\item Every isomorphism \(\theta:\M_{\mathcal U,\mathbf d}\to\M_{\mathcal V,\mathbf e}\) sends \(D\) to a
unitary conjugate of \(F\).
\item If \(\M_{\mathcal U,\mathbf d}\equiv\M_{\mathcal V,\mathbf e}\), then \((\M_{\mathcal U,\mathbf d},D)\equiv(\M_{\mathcal V,\mathbf e},F)\).
\item If \(\M_{\mathcal U,\mathbf d}\equiv\M_{\mathcal V,\mathbf e}\), then
\[
 \prod_{\mathcal U}^{\mathrm{met}}\operatorname{Sym}(d_n)
 \ \equiv\
 \prod_{\mathcal V}^{\mathrm{met}}\operatorname{Sym}(e_n),
\]
where the symmetric groups carry the normalized Hamming metric.
\end{enumerate}
\end{maintheorem}

To obtain the last assertion, permutation matrices identify the
Weyl group with the metric ultraproduct of
symmetric groups. Its Hamming metric is one half the square of
the quotient \(2\)-metric. Definability lets us make this
identification under elementary equivalence, and
\cite[Theorem~5.6 and Remark~5.7]{AT} then gives the matrix-algebra
assertion of Theorem~\ref{cor:many-theories}. Reconstruction from
unitary groups gives the other assertion. The argument that preserves the
diagonal also extends to arbitrary internal algebras, giving the
following structural conclusion.

\begin{maintheorem}\label{thm:invariance}
Let \(\theta:\M_{\mathcal U,\mathbf d}\to\M_{\mathcal V,\mathbf e}\)
be an isomorphism of tracial matrix ultraproducts. Then for every
unital von Neumann subalgebra
\(A\subseteq\M_{\mathcal U,\mathbf d}\),
\[
 A\text{ is internal}
 \quad\Longleftrightarrow\quad
 \theta(A)\text{ is internal}.
\]
\end{maintheorem}

The constructions used here also give fixed finite distance formulas
over four unitary parameters for internal abelian algebras and six
for arbitrary internal algebras; see
Proposition~\ref{thm:internal-finite-parameter-definability}.
They also show that, for commuting internal algebras \(A,B\), both
\(A'\cap B'\) and their relative join \((A'\cap B')'\) are internal
(Corollary~\ref{prop:gen-commuting-internal}).

We next use the recovered permutation groups to distinguish an
explicit family of factors from their amplifications. For a \(\mathrm{II}_1\)
factor \(Q\), let \(Q^t=pM_k(Q)p\), with normalized trace,
where \(k>t\) and \((\tau_k\otimes\tau)(p)=t/k\). Its isomorphism
class depends only on \(t>0\). The \emph{first-order fundamental group}
\[
 \mathcal F_{\mathrm{fo}}(Q)=\{t>0:Q^t\equiv Q\}
\]
contains \(\mathcal F(Q)=\{t>0:Q^t\cong Q\}\); see
\cite[Section~4]{GoldbringHartSurvey}.
For the construction, put \(G_r=\mathrm{PSL}_2(2^r)\), for \(r\geq2\), so that
\(|G_r|=2^r(2^{2r}-1)\).

\begin{maintheorem}\label{thm:fundamental-group}
Fix a nonprincipal ultrafilter \(\mathcal U\) on \(\N\).
For each set \(S\) of odd primes, put
\[
 r_n(S)=2^n\prod_{\substack{p\in S\\p\leq n}}p,
 \qquad
 Q_S=\prod_{\mathcal U}(M_{|G_{r_n(S)}|}(\C),\tau_{|G_{r_n(S)}|}).
\]
Then \(\mathcal F_{\mathrm{fo}}(Q_S)=\mathcal F(Q_S)=\{1\}\).
Distinct sets \(S\) give non-elementarily equivalent factors.
\end{maintheorem}

To prove this theorem, we follow the left-right action method of
\cite[Section~5]{AT} and recover the regular actions of
\(\mathrm{PSL}_2(2^r)\) from their elementary properties.
The recovered actions force the ambient dimension to be asymptotic
to a group order. We then use element orders to detect divisibility
of \(r\), which rules out nontrivial amplifications.


Popa asked in 2008 whether there are uncountably many nonisomorphic
tracial matrix ultraproducts \cite[Introduction]{FHS1}.
For the related ultrapower problem, Ge and Hadwin proved that,
under CH, all tracial ultrapowers of a fixed finite von Neumann
algebra with separable predual are isomorphic
\cite[Theorem~3.2]{GeHadwin}. Li later suggested that
matrix ultraproducts might be independent of the ultrafilter
\cite[p.~1287]{Li}. Under failure of CH, Farah, Hart and Sherman
obtained nonisomorphic matrix ultraproducts
\cite[Proposition~3.3]{FHS1}, and Farah and Shelah obtained the
maximal number \(2^{\mathfrak c}\) \cite[Proposition~8.3]{FarahShelah}.
Under CH, matrix ultraproducts are isomorphic exactly when they
are elementarily equivalent. Theorem~\ref{cor:many-theories}
distinguishes their theories in ZFC and hence answers Popa's
question and gives continuum many classes
in the classification problem of \cite[Problem~1]{Kunnawalkam}.
For the model-theoretic framework and the game formulation, see
\cite{FHS2,FHS3} and \cite[Question~4.1]{FarahGames}.

For universal sofic groups, Thomas obtained the maximal number of
isomorphism types under failure of CH and realized algebraic
ultraproducts of finite expander groups as centralizers
\cite[Theorems~1.1 and~2.2]{Thomas}. The authors obtained continuum
many theories in ZFC in \cite{AT} and proposed the diagonal-normalizer
route to matrix ultraproducts \cite[Remark~5.8]{AT}. Jekel discussed
this route in \cite[Section~1.1]{JekelInformation}.
Since full types include sentence values, nonconvergence also
obstructs convergence of full types of random matrix tuples.
It answers the corresponding ultrafilter-independence question for
noncommutative stochastic processes as well; see
\cite[Question~5.6]{JekelSurvey}, \cite[Section~1.1]{JekelInformation},
and \cite[Remark~1.3 and Proposition~D(1)]{JekelProcesses}.

Our internality arguments use Popa's coordinate commutant identities
\cite[Section~1.6]{PopaIndependence} and conjugacy of coordinate
MASAs \cite[Section~2.1(a)]{PopaKS}. The diagonal is a MASA,
but its normalizer generates a proper subalgebra: tracial matrix
ultraproducts have no Cartan subalgebras
\cite[Proposition~2.3]{PopaKS}; see also
\cite[Theorem~7.3]{PopaOrthogonal} for ultrapowers of factors.

\medskip

While preparing this manuscript, we learned of the independent work
of Ando and Goldbring \cite{AndoGoldbring}. They prove that every
tracial matrix ultraproduct has trivial first-order fundamental group,
which also implies nonconvergence.

\medskip

The paper is organized as follows: Section~\ref{sec:preliminaries} supplies the internality results
needed to recognize the diagonal. Using these,
Section~\ref{sec:defining-diagonal} recovers the diagonal and its
normalizer under isomorphism, and
Section~\ref{sec:internal-isomorphisms} extends preservation to
all internal algebras. To pass to elementary equivalence, we make
the constructions definable in Section~\ref{sec:logical-details}.
Section~\ref{sec:left-right} then recovers the left-right actions
from admissible tuples and uses element orders to prove
Theorem~\ref{thm:fundamental-group}.

\section{Spectral gap and internal algebras}\label{sec:preliminaries}

We use the normalized trace \(\tau_d=d^{-1}\operatorname{Tr}\) on
\(M_d=M_d(\mathbb C)\) and its norm
\(\|x\|_{2,d}=\tau_d(x^*x)^{1/2}\). We write
\(A_{\leq1}\) for the operator-norm unit ball of an algebra \(A\).
Bounded functional calculus lets us choose contractive coordinate
representatives for every unit-ball element of an internal algebra.

We use primes for commutants in the specified ambient matrix
ultraproduct. For a finite tuple \(\mathbf u\), we write \(\{\mathbf u\}'\cap\M_{\mathcal U,\mathbf d}\) for its commutant, and we define the relative
join of two subalgebras by
\begin{equation}\label{eq:join}
 A\vee_{\mathrm{rel}}B=(A'\cap B')'.
\end{equation}
Thus the relative join is the smallest algebra containing \(A\)
and \(B\) that equals its relative bicommutant.

\subsection{Internal subalgebras, commutants and spectral gap}
\label{sec:internal-preliminaries}
The following identities let us pass between coordinate models and
commutants in the ultraproduct. They also reduce unitary conjugacy
of internal MASAs to the matrix case.

\begin{lemma}\label{lem:commutant}
Let \(\M_{\mathcal U,\mathbf d}\) be a tracial matrix ultraproduct,
and let \(A_n\subseteq M_{d_n}(\mathbb C)\) be unital \(*\)-subalgebras.
Put \(A=[A_n]_{\mathcal U}\). Then
\begin{align*}
 A'&=[A_n'\cap M_{d_n}]_\mathcal U, \qquad
 Z(A)=[Z(A_n)]_\mathcal U,\qquad
 A''=A.
\end{align*}
The coordinate conditional expectations \(E_{A_n}:M_{d_n}(\mathbb C)\to A_n\)
induce the conditional expectation \(E_A:\M_{\mathcal U,\mathbf d}\to A\).
Every internal abelian algebra has abelian coordinate models and is
contained in an internal MASA. Every internal MASA is a MASA in the ambient ultraproduct and is
unitarily conjugate to the canonical diagonal \([D_{d_n}]_\mathcal U\).
\end{lemma}

\begin{proof}
We first identify the conditional expectation induced by the coordinate
maps. Since each \(E_{A_n}\) is a contraction in both the operator norm
and the \(2\)-norm, these maps induce a trace-preserving projection
onto \([A_n]_\U\). The trace-pairing characterization of conditional
expectation then identifies this projection with \(E_A\).

The commutant identity \(A'=[A_n'\cap M_{d_n}]_\U\) follows from \cite[Section~1.6]{PopaIndependence}. To obtain the
center identity, we apply the same result with ambient coordinate
algebras \(A_n\), which gives \(Z(A)=[Z(A_n)]_\U\). Applying the commutant identity
twice in \(M_{d_n}\) and using the finite-dimensional bicommutant theorem
then gives \(A''=A\).

Suppose now that \(A\) is abelian. The center identity allows us
to replace its coordinate models by \(Z(A_n)\). We can therefore extend
each model to a MASA \(D_n^0\subseteq M_{d_n}\), and the commutant identity shows that
\([D_n^0]_\U\) is a MASA containing \(A\). If \(A\) is maximal among
internal abelian algebras, this inclusion forces \(A=[D_n^0]_\U\), so it is
also maximal abelian in the ambient ultraproduct. Conversely, an
internal algebra that is maximal abelian is maximal among internal
abelian algebras. Finally, coordinate unitary conjugacies identify
each of these MASAs with the canonical diagonal
\cite[Section~2.1(a)]{PopaKS}. Its maximal abelianness can also be read directly from
the commutant identity:
\[
 [D_{d_n}]_\U'\cap\M_{\mathcal U,\mathbf d}
 =[D_{d_n}'\cap M_{d_n}]_\U=[D_{d_n}]_\U.
\]
\end{proof}

To recognize an algebra given as a commutant, we will use spectral
gap. Let \(\mathcal P\) be a tracial von Neumann algebra and let
\(\mathbf u=(u_1,\ldots,u_h)\in U(\mathcal P)^h\), with \(h\geq1\).
Put \(\mathcal C=\mathcal P\cap\{u_1,\ldots,u_h\}'\). We say that
\(\mathbf u\) has \emph{spectral gap in \(\mathcal P\)} if there is
\(0<L<\infty\) such that
\begin{equation}\label{eq:sgi}
 \|x-E_{\mathcal C}x\|_2^2
 \leq L\sum_{j=1}^h\|[u_j,x]\|_2^2
 \qquad(x\in\mathcal P).
\end{equation}
We call any such \(L\) a \emph{spectral gap constant} for the tuple.
The distance in this estimate is measured to the commutant within
the specified ambient algebra.

Our main input, the spectral gap internality theorem
\cite[Theorem~B]{ALT}, turns this estimate into coordinate models
for the commutant. We use its relative form so that we can also
work inside an internal subalgebra.

\begin{theorem}
\label{thm:spectral-gap-internality}
Let \(\mathcal P=[\mathcal P_n]_\U\) be an internal subalgebra of a
tracial matrix ultraproduct \(\M_{\mathcal U,\mathbf d}\), and let
\(u_1,\ldots,u_h\in U(\mathcal P)\) have spectral
gap in \(\mathcal P\). Then their relative commutant
\(\mathcal C=\mathcal P\cap\{u_1,\ldots,u_h\}'\) has the form
\(\mathcal C=[\mathcal C_n]_\U\) for unital \(*\)-subalgebras
\(\mathcal C_n\subseteq\mathcal P_n\).
\end{theorem}

\begin{proof}
We will approximate the conditional expectation onto the relative
commutant by expectations on the coordinate summands, where the
matrix version of the theorem applies. Let \(L\) be a spectral gap
constant, and decompose the coordinate algebras with their inherited
traces as
\[
 (\mathcal P_n,\tau_{d_n}|_{\mathcal P_n})
 =\bigoplus_s(M_{a_{n,s}},w_{n,s}\tau_{a_{n,s}}),
 \qquad w_{n,s}>0,\quad\sum_s w_{n,s}=1.
\]
We choose unitary lifts \(u_{j,n}\in\mathcal P_n\) and form their lazy conjugation averages
\[
 T_n(x)=\frac12x+\frac1{4h}\sum_{j=1}^h
       (u_{j,n}xu_{j,n}^*+u_{j,n}^*xu_{j,n}).
\]
The diagonal argument of \cite[Proposition~2.2]{ALT} applies to finite-dimensional tracial
algebras. It therefore gives powers \(F_n=T_n^{k_n}\), with \(k_n\to_\U\infty\), that induce
\(E_{\mathcal C}\). Since these maps preserve every summand, we may work with
their restrictions \(T_{n,s},F_{n,s}\). We use the corresponding restrictions
for the unitary entries as well.

To find suitable expectations on the summands, fix \(\varepsilon>0\) and call
a summand bad if \(F_{n,s}\) has \(\infty\)-to-\(2\) distance at least
\(\varepsilon\) from every expectation onto a unital subalgebra of \(M_{a_{n,s}}\).
We claim that the total weight of the bad summands tends to zero.
To measure the failure of the spectral gap estimate on each summand,
we set
\[
 \eta_{n,s}=\sup_{\|x\|\leq1}
 \left(\|x-F_{n,s}x\|_2^2
       -L\sum_{j=1}^h\|[u_{j,n,s},x]\|_2^2\right)\geq0,
\]
where we use the normalized trace on that summand. Choosing maximizing
contractions in all the summands and applying spectral gap in
\(\mathcal P\) gives \(\sum_s w_{n,s}\eta_{n,s}\to_\U0\). If the total bad weight did not tend to zero,
Markov's inequality would then allow us to choose bad summands
\(s(n)\) for which \(\eta_{n,s(n)}\to_\U0\).

We show that this choice contradicts the matrix version of the theorem.
In the auxiliary tracial ultraproduct \(Q=\prod_\U M_{a_{n,s(n)}}\), put \(v_j=[u_{j,n,s(n)}]_\U\) and
\(D=Q\cap\{v_1,\ldots,v_h\}'\). The errors on the selected summands tend to zero, so the
induced maps \(T,F\) satisfy
\[
 \|x-Fx\|_2^2\leq L\sum_{j=1}^h\|[v_j,x]\|_2^2.
\]
We also have \(TF=F\): indeed, \(0\leq T_{n,s}\leq I\) on \(L^2\) gives the bound
\(\|(I-T_{n,s})T_{n,s}^{k_n}\|_{2\to2}\leq1/(k_n+1)\). It follows that \(F(Q)\subseteq D\), while the displayed inequality
shows that \(F\) fixes \(D\) pointwise. Since the induced map
is self-adjoint, we obtain \(F=E_D\). Thus the same inequality is a
spectral gap estimate, and we can apply \cite[Theorem~B]{ALT}, whose statement
also allows bounded coordinate dimensions, to obtain \(D=[D_n]_\U\).
The coordinate expectations and the selected powers now induce the
same map. Hence \cite[Lemma~2.1]{ALT} gives \(\|F_{n,s(n)}-E_{D_n}\|_{\infty\to2}\to_\U0\), contrary to the choice of
bad summands.

We have proved that the bad weight tends to zero for every \(\varepsilon>0\).
In each summand, choose an expectation within
\(1/n\) of the best approximation to \(F_{n,s}\), and let \(\mathcal C_n\subseteq\mathcal P_n\) be the
direct sum of the resulting ranges. All errors are bounded by \(2\),
so separating the good and bad summands gives
\[
 \|F_n-E_{\mathcal C_n}\|_{\infty\to2}^2
 \leq(\varepsilon+1/n)^2+4\sum_{s\ \mathrm{bad}}w_{n,s}.
\]
Taking the ultralimit and then letting \(\varepsilon\) tend to zero shows
that the total error vanishes. Thus these coordinate expectations
induce \(E_{\mathcal C}\), and we conclude that \(\mathcal C=[\mathcal C_n]_\U\).
\end{proof}

\subsection{Commuting internal abelian algebras}\label{sec:commuting-internal}
To combine commuting internal abelian algebras, we need their
relative join to remain internal. Their coordinate models need
not commute exactly, so we use the abelian stability theorem
\cite[Theorem 1.6]{AD}. It states
that for every \(\delta>0\) there is \(\varepsilon>0\), independent of the
countable abelian group \(G\) and the dimension \(d\), such that a map
\(f:G\to U(d)\) with uniform multiplicative defect less than \(\varepsilon\) in
normalized Hilbert--Schmidt norm is uniformly within \(\delta\) of a unitary
representation on \(\C^d\).

\begin{proposition}\label{prop:abelian-commuting-internal}
Let \(A,B\) be commuting internal abelian subalgebras of a tracial
matrix ultraproduct \(\M_{\mathcal U,\mathbf d}\). Then \(A'\cap B'\)
and \(A\vee_{\mathrm{rel}}B\) are internal, and the latter is abelian.
\end{proposition}
\begin{proof}
We choose abelian models \(A=[A_n]_{\U}\) and \(B=[B_n]_{\U}\) and first record what
commutation means in the coordinates. Since the two internal algebras
commute, a coordinate selection argument gives
\[
 \varepsilon_n:=\sup_{a\in U(A_n),\ b\in U(B_n)}
 \|[a,b]\|_{2,d_n}\longrightarrow_{\U}0.
\]
To replace this approximate commutation by an exact abelian action,
we choose countable dense subgroups \(G_n\subseteq U(A_n)\) and \(H_n\subseteq U(B_n)\). The map
\(f_n:G_n\times H_n\to U(d_n)\) defined by \(f_n(a,b)=ab\) then has multiplicative defect at most
\(\varepsilon_n\), since its defect on \((a,b),(c,e)\) is \(a[b,c]e\). Uniform abelian
stability therefore gives representations \(\pi_n\) in the same
dimensions satisfying
\[
 \eta_n:=\sup_{a\in G_n,\ b\in H_n}\|\pi_n(a,b)-ab\|_{2,d_n}
 \longrightarrow_{\U}0.
\]
We use their commutants as coordinate models by setting \(C_n=\pi_n(G_n\times H_n)'\cap M_{d_n}\).
To identify the resulting ultraproduct, take \(x=[x_n]_{\U}\in A'\cap
B'\) with bounded
representatives. A coordinate selection argument shows that these
representatives commute asymptotically and uniformly with \(G_n\)
and \(H_n\). By the uniform approximation just obtained, they also
commute asymptotically and uniformly with \(\pi_n(G_n\times
H_n)\). Averaging over
the compact closure of this image with respect to Haar measure
therefore gives
\[
 \|x_n-E_{C_n}x_n\|_{2,d_n}
 \leq\sup_{g\in G_n\times H_n}\|[x_n,\pi_n(g)]\|_{2,d_n}
 \longrightarrow_{\U}0.
\]
Thus these representatives can be replaced by elements of the
coordinate commutants without changing their ultraproduct class.

Conversely, representatives in \(C_n\) commute asymptotically and
uniformly with \(G_n\) and \(H_n\) by comparison with \(\pi_n(a,1)\) and
\(\pi_n(1,b)\). Density then shows that their ultraproduct classes commute
with \(A\) and \(B\). We have proved \(A'\cap B'=[C_n]_{\U}\), and Lemma~\ref{lem:commutant}
shows that its commutant is internal as well. To check that this
relative join is abelian, put \(S=A\cup B\). Since the original algebras
are abelian and commute, we have \(S\subseteq S'\), and taking commutants gives
\(S''\subseteq S'=(S'')'\). This is precisely the required abelianness.
\end{proof}

\begin{corollary}\label{lem:envelope}
Every finite family of pairwise commuting internal abelian subalgebras
of a tracial matrix ultraproduct \(\M_{\mathcal U,\mathbf d}\) is
contained in a common internal MASA.
\end{corollary}
\begin{proof}
We apply Proposition~\ref{prop:abelian-commuting-internal} successively to form the relative joins
of the given algebras. At each step, the next algebra commutes with
the relative bicommutant of the preceding family, so the proposition
applies again and the resulting join remains internal and abelian.
The final join contains the whole family, and Lemma~\ref{lem:commutant} places
it in an internal MASA.
\end{proof}

\section{Recovering the diagonal and its normalizer}\label{sec:defining-diagonal}
\label{sec:short-consequences}

We construct two tuples with spectral gap whose abelian commutants
have the canonical diagonal as their joint commutant. The internality
results then identify its image under an isomorphism. Passing to
the Weyl group recovers the associated universal sofic
group and gives the isomorphism conclusions.

\subsection{Explicit quantum expanders}\label{sec:explicit-expanders}

We need unitaries whose spectral gap constants do not depend on
dimension. The following explicit pairs give such estimates for
scalar commutants and for sums of scalar blocks of distinct sizes.

\begin{lemma}\label{lem:expander-gaps}
For \(k\geq1\), let \(\rho_k=\operatorname{Sym}^{k-1}(\C^2)\)
be the \(k\)-dimensional irreducible unitary representation of \(\mathrm{SU}(2)\),
written in the normalized symmetric-tensor basis, and put
\begin{gather*}
 g_1=\frac1{\sqrt3}\begin{pmatrix}1+i&1\\-1&1-i\end{pmatrix},\qquad
 g_2=\frac1{\sqrt3}\begin{pmatrix}1+i&-1\\1&1-i\end{pmatrix},\\
 v_{j,k}=\rho_k(g_j)\quad(j=1,2).
\end{gather*}
Then the following hold.
\begin{enumerate}
\item For every \(k\geq1\), the pair \((v_{1,k},v_{2,k})\) has
commutant \(\C1\) in \(M_k\) and spectral gap with constant \(2\):
\[
 \|z-\tau_k(z)1\|_{2,k}^2
 \leq2\sum_{j=1}^2\|[z,v_{j,k}]\|_{2,k}^2
 \qquad(z\in M_k).
\]
\item Let \(m_1,\ldots,m_r\) be pairwise distinct positive integers, and put
\begin{align*}
 C&=\bigoplus_{i=1}^r\C1_{m_i}\subseteq M_d,
 \qquad d=\sum_i m_i,\\
 \widehat v_j&=\bigoplus_{i=1}^r v_{j,m_i}\quad(j=1,2).
\end{align*}
Then \((\widehat v_1,\widehat v_2)\) has commutant \(C\) in \(M_d\)
and spectral gap with constant \(8\):
\[
 \|x-E_Cx\|_2^2\leq8\sum_{j=1}^2\|[x,\widehat v_j]\|_2^2
 \qquad(x\in M_d).
\]
\end{enumerate}
\end{lemma}

\begin{proof}
To prove (1), we relate the given unitaries to the Ramanujan
construction. Under the usual identification of unit quaternions
with \(\mathrm{SU}(2)\), they correspond to
\(g_1=(1+\mathbf i+\mathbf j)/\sqrt3\) and
\(g_2=(1+\mathbf i-\mathbf j)/\sqrt3\). For this comparison, put
\(s_{\varepsilon,\eta}=(\mathbf k+\varepsilon\mathbf i+
\eta\mathbf j)/\sqrt3\), for \(\varepsilon,\eta\in\{-1,1\}\).
These are the four half-turns in the \(p=3\) construction of
\cite[Appendix~C]{IJJS}. Since
\[
 \{\mathbf k s_{\varepsilon,\eta}:\varepsilon,\eta\in\{-1,1\}\}
 =\{-g_1,-g_1^{-1},-g_2,-g_2^{-1}\},
\]
we can apply the Ramanujan estimate in
\cite[Section~6.1, (93)]{IJJS} to obtain
\begin{align*}
 \left\|T_k\big|_{\ker\tau_k}\right\|_{2\to2}&\leq2\sqrt3, \qquad
 T_k=\sum_{j=1}^2\bigl(\operatorname{Ad}(v_{j,k})+
                         \operatorname{Ad}(v_{j,k}^*)\bigr).
\end{align*}
Indeed, multiplication by \(\operatorname{Ad}(\rho_k(\mathbf k))\)
preserves the norm, while the central signs disappear under
conjugation. To obtain the required commutator estimate, take
\(z\in M_k\) and \(z_0=z-\tau_k(z)1\). The preceding bound yields
\[
 \sum_{j=1}^2\|[z,v_{j,k}]\|_{2,k}^2
 =\langle z_0,(4I-T_k)z_0\rangle
 \geq(4-2\sqrt3)\|z_0\|_{2,k}^2\geq\tfrac12\|z_0\|_{2,k}^2.
\]
This gives the spectral gap estimate with constant \(2\). It also shows that every element of the commutant
lies in \(\C1\), completing the proof of (1).

For (2), we first use (1) to obtain a gap for arbitrary
representations. Given a unit vector \(\xi\) in a nontrivial
irreducible representation \(\rho_m\), we apply the estimate to
the rank-one projection \(p\) onto \(\C\xi\).
Writing \(\|x\|_{\mathrm{HS}}=\operatorname{Tr}(x^*x)^{1/2}\)
for the unnormalized Hilbert--Schmidt norm, we obtain
\[
 \frac12\leq1-\frac1m
 =\left\|p-\frac1m1\right\|_{\mathrm{HS}}^2
 \leq2\sum_{j=1}^2\|[p,\rho_m(g_j)]\|_{\mathrm{HS}}^2
 \leq4\sum_{j=1}^2\|\rho_m(g_j)\xi-\xi\|^2.
\]
Here the first inequality uses \(m\geq2\), and the last follows
from the rank-one calculation
\[
 \|[p,v]\|_{\mathrm{HS}}^2
 =2\bigl(1-|\langle v\xi,\xi\rangle|^2\bigr)
 \leq2\|v\xi-\xi\|^2
\]
for a unitary \(v\). The resulting bound is independent of the
irreducible representation. By homogeneity and complete reducibility,
it follows that every finite-dimensional continuous unitary
representation \(\pi\) of \(\mathrm{SU}(2)\) satisfies
\begin{equation}\label{eq:direct-representation-gap}
 \|\xi-p_\pi\xi\|^2
 \leq8\sum_{j=1}^2\|\pi(g_j)\xi-\xi\|^2,
\end{equation}
where \(p_\pi\) projects onto the invariant vectors.

We now apply \eqref{eq:direct-representation-gap} to the conjugation
representation associated with \(\rho=\bigoplus_i\rho_{m_i}\) on
\(L^2(M_d)\). The summands \(\rho_{m_i}\) are pairwise inequivalent,
so Schur's lemma identifies the commutant as
\(\rho(\mathrm{SU}(2))'=C\). Thus averaging the conjugation action
gives \(E_C\), and the representation gap becomes precisely the
estimate in (2). Since vanishing of its right-hand side forces
membership in this algebra, it also identifies the commutant of
\(\widehat v_1,\widehat v_2\) with \(C\).
\end{proof}

\subsection{Recovering the diagonal}

The scalar estimate of Lemma~\ref{lem:expander-gaps} lets us obtain
row and column algebras as commutants with spectral gap. Their
joint commutant is the diagonal. The block estimate will be needed
for general internal algebras in Section~\ref{sec:internal-isomorphisms}.

\begin{lemma}\label{lem:row-column-recovery}
For every \(k\geq1\), there are two unitary triples
\((z_{i1,k},z_{i2,k},z_{i3,k})\), \(i=1,2\), in
\(M_{k^2}=M_k\otimes M_k\), with commutants
\[
 B_{1,k}=D_k\otimes1,\qquad B_{2,k}=1\otimes D_k,
\]
and spectral gap with constant \(10\):
\[
 \|x-E_{B_{i,k}}x\|_2^2
 \leq10\sum_{j=1}^3\|[x,z_{ij,k}]\|_2^2
 \qquad(x\in M_{k^2},\ i=1,2).
\]
For every tracial matrix ultraproduct \(\M_{\mathcal U,\mathbf d}\),
these triples yield six unitaries \(z_{ij}\) whose commutants
\(B_i=\{z_{i1},z_{i2},z_{i3}\}'\) are commuting internal
abelian algebras, have the same spectral gap estimates, and satisfy
\[
 B_1'\cap B_2'=[D_{d_n}]_\U.
\]
\end{lemma}

\begin{proof}
We construct the first triple by using the expander pair to leave
one matrix factor in the commutant, and a third unitary to reduce
that factor to its diagonal. Index the standard basis of \(\C^k\)
by \(0,\ldots,k-1\), let \(S\) be its cyclic shift, and set
\[
 W=\sum_{r=0}^{k-1}e_{rr}\otimes S^r,\qquad
 (z_{11,k},z_{12,k},z_{13,k})=(1\otimes v_{1,k},1\otimes v_{2,k},W).
\]
We take the tensor flip for the second triple, so it suffices to
verify the first. Put \(Q=M_k\otimes1\) and
\(a=\|x-E_Qx\|_2\). Tensoring Lemma~\ref{lem:expander-gaps}(1)
with the identity gives
\(a^2\leq2\sum_{j=1}^2\|[x,z_{1j,k}]\|_2^2\), which controls
the distance to the first matrix factor. Within that factor, the
orthonormality of the powers of \(S\) for the normalized trace gives
\[
 \|[y,W]\|_2^2=2\|y-E_{B_{1,k}}y\|_2^2\qquad(y\in Q).
\]
We combine these two estimates by writing \(c=\|[x,W]\|_2\)
and applying Pythagoras:
\[
 \|x-E_{B_{1,k}}x\|_2^2
 \leq a^2+\tfrac12(c+2a)^2
 \leq5a^2+c^2
 \leq10\sum_{j=1}^3\|[x,z_{1j,k}]\|_2^2.
\]
This proves the asserted gap. Since the triple commutes with
\(B_{1,k}\), the same estimate identifies its commutant, and the
tensor flip gives the corresponding conclusions for the second triple.

To pass to arbitrary \(d_n\), we put
\(k_n=\lfloor\sqrt{d_n}\rfloor\) and extend the six unitaries by
the identity on the complement of the standard
\(k_n^2\)-dimensional corner \(p_n\). The relative dimension of
this corner tends to one, so
\[
 \|x_n-p_nx_np_n\|_{2,d_n}
 \leq2\|x_n\|\sqrt{1-k_n^2/d_n}\longrightarrow_\U0
\]
for every bounded sequence. Hence corner inclusion identifies the
ultraproduct of the square matrices with
\(\M_{\mathcal U,\mathbf d}\). This identification preserves the
canonical diagonals and transfers the gap estimates. In particular,
the resulting \(B_i\) are the ultraproducts of the row and column
algebras, extended by scalars on the negligible complement.

It remains to identify their joint commutant. In each square corner,
the expectations onto \(B_{1,k}'=D_k\otimes M_k\) and
\(B_{2,k}'=M_k\otimes D_k\) commute and have product
\(E_{D_{k^2}}\). Their product estimate therefore passes to the
ultraproduct, where Lemma~\ref{lem:commutant} identifies the two
commutants and gives the asserted intersection.
\end{proof}

The next lemma recognizes an internal diagonal from these spectral
gap and commutation conditions. This will let us transport the
construction under isomorphisms and, later, by transferring its
defining conditions.

\begin{lemma}\label{lem:diagonal-recognition}
Let \(C_1,C_2\) be commuting abelian subalgebras of a tracial
matrix ultraproduct. Suppose that each \(C_i\) is the commutant
of a finite unitary tuple with spectral gap, and that
\(D=C_1'\cap C_2'\) is abelian. Then \(D\) is an internal
MASA, and hence is unitarily conjugate to the canonical diagonal.
\end{lemma}

\begin{proof}
The spectral gap hypothesis makes both \(C_i\) internal by
Theorem~\ref{thm:spectral-gap-internality}. We first identify their
joint commutant as a MASA: since the algebras commute and are abelian,
\(C_i\subseteq D\), whence \(D'\subseteq D\). Together with the
assumed abelianness, this shows that \(D\) is a MASA.
To prove its internality, apply Corollary~\ref{lem:envelope} to
place \(C_1,C_2\) in a common internal MASA \(E\).
Then \(E\subseteq D\), and maximal abelianness forces \(E=D\).
The required conjugacy now follows from Lemma~\ref{lem:commutant}.
\end{proof}

\begin{proof}[Proof of Theorem~\ref{thm:reduction}(1)]
Let \(\theta:\M_{\mathcal U,\mathbf d}\to\M_{\mathcal V,\mathbf e}\)
be an isomorphism, and let \(D,F\) be the canonical diagonals.
To apply the recognition lemma to the image diagonal, we use
Lemma~\ref{lem:row-column-recovery} to write \(D=B_1'\cap B_2'\).
Since the isomorphism preserves the trace and commutants, it also
preserves the spectral gap estimates. Consequently
\(\theta(B_1),\theta(B_2)\) satisfy the hypotheses of
Lemma~\ref{lem:diagonal-recognition}, with joint commutant
\(\theta(D)\). We conclude from that lemma that \(\theta(D)\)
is an internal MASA and hence is unitarily conjugate to \(F\).
\end{proof}

\subsection{Weyl groups and isomorphisms}\label{sec:isomorphism-consequences}

Having recovered the diagonal up to conjugacy, we use its normalizer
to recover the universal sofic group. For an internal MASA
\(D\subseteq\M_{\mathcal U,\mathbf d}\), write
\begin{align*}
 N_{\M_{\mathcal U,\mathbf d}}(D)
 &=\{u\in U(\M_{\mathcal U,\mathbf d}):uDu^*=D\},\\
 \mathcal W(\M_{\mathcal U,\mathbf d},D)
 &=N_{\M_{\mathcal U,\mathbf d}}(D)/U(D).
\end{align*}
We call \(\mathcal W(\M_{\mathcal U,\mathbf d},D)\) the \emph{Weyl group}.
We identify permutations with their permutation matrices.
In the coordinate algebras, every element of \(N_{M_n}(D_n)\) has the form \(\sigma w\),
where \(\sigma\in\Sym(n)\) and \(w\in U(D_n)\).
To carry this description to the ultraproduct, we need to represent
each normalizer by coordinate normalizers. The next estimate gives
the required approximation.

Popa, Sinclair and Smith proved quantitative estimates for the
distance of a unitary to the normalizer of a MASA in a separably
acting \(\mathrm{II}_1\) factor \cite[Theorem~6.2]{PopaSinclairSmith}.
For matrix diagonals, the following elementary argument gives the
bound needed here.

\begin{lemma}\label{lem:normalizer}
For \(u\in U(n)\), put
\[
 r_{D_n}(u)=\sup_{a\in(D_n)_{\leq1}}
       \|uau^*-E_{D_n}(uau^*)\|_2.
\]
There is a unitary \(v\in N_{M_n}(D_n)\) such that
\begin{equation}\label{eq:normalizer}
 \|u-v\|_2\leq\sqrt2\,r_{D_n}(u).
\end{equation}
The same estimate holds for an internal MASA in a tracial matrix
ultraproduct \(\M_{\mathcal U,\mathbf d}\). After conjugating it to
the canonical diagonal, each unitary normalizer has coordinate representatives that normalize
the coordinate diagonals.
\end{lemma}

\begin{proof}
We use the doubly stochastic matrix \(P=(|u_{ij}|^2)_{i,j}\) to
choose the support of a nearby normalizing unitary. To relate its
entries to the normalizer defect, take independent uniform signs
\(\varepsilon_j\) and set
\(s=\operatorname{diag}(\varepsilon_1,\ldots,\varepsilon_n)\).
Averaging over these signs gives
\[
 \mathbb E\|usu^*-E_{D_n}(usu^*)\|_2^2
 =1-\frac1n\sum_{i,j}P_{ij}^2.
\]
Thus, with \(r=r_{D_n}(u)\), we obtain
\(\sum_{i,j}P_{ij}^2\geq n(1-r^2)\). By Birkhoff's theorem,
\(P\) is a convex combination of permutation matrices, so one
of the corresponding permutations, say \(\sigma\), satisfies
\[
 \sum_iP_{i,\sigma(i)}\geq\sum_{i,j}P_{ij}^2.
\]
We choose \(v\in N_{M_n}(D_n)\) supported on these positions
and give its nonzero entries the phases of the corresponding
entries of \(u\). At a zero entry of \(u\), the phase can be
chosen arbitrarily. With this choice, the preceding bound yields
\[
 \|u-v\|_2^2
 =2-\frac2n\sum_i|u_{i,\sigma(i)}|
 \leq2-\frac2n\sum_iP_{i,\sigma(i)}
 \leq2r^2.
\]
This proves the matrix estimate.

We next pass the estimate to a tracial matrix ultraproduct. Polar
decomposition allows us to lift its unitaries to unitary sequences.
For \(D=[D_{d_n}]_\U\), the coordinate expectations induce
\(E_D\), and the normalizer defects satisfy
\[
 r_D([u_n]_\U)=\lim_{\U}r_{D_{d_n}}(u_n).
\]
Indeed, the coordinate suprema give one inequality, and choosing a
maximizer in each coordinate gives the reverse inequality. The
matrix estimate therefore passes to the ultraproduct. If \(u\)
normalizes \(D\), the defect vanishes, so the resulting coordinate
diagonal normalizers represent \(u\). Finally, unitary conjugacy
transfers both conclusions to every internal MASA.
\end{proof}

For the canonical diagonal, Lemma~\ref{lem:normalizer} now lets us
identify the Weyl group with the associated universal
sofic group:
\begin{equation}\label{eq:weyl}
 \mathcal W(\M_{\mathcal U,\mathbf d},[D_{d_n}]_\U)
 \cong\prod_\U^{\mathrm{met}}\Sym(d_n).
\end{equation}
Indeed, every normalizer has coordinate representatives
\(\sigma_n w_n\), and the diagonal factors disappear in the
quotient. To identify the metric as well, we observe that for
permutation matrices,
\begin{equation}\label{eq:quotmetric}
 \frac12\inf_{w\in U(D_n)}\|\sigma-\rho w\|_2^2
 =1-\frac{\#\{j:\sigma(j)=\rho(j)\}}n
 =d(\sigma,\rho).
\end{equation}
Indeed, we match phases at the common nonzero entries, while the
other entries have disjoint supports. We equip \(\mathcal W\) with the metric
\[
 d(\bar u,\bar v)=\frac12\inf_{w\in U(D)}\|u-vw\|_2^2.
\]
Since coordinate minima pass to ultralimits, this is the normalized
Hamming metric under the identification above.

\begin{isomorphismtheorem}\label{thm:isomorphism-reduction}
Let \(D=[D_{d_n}]_\U\subseteq\M_{\mathcal U,\mathbf d}\) and
\(F=[D_{e_n}]_\V\subseteq\M_{\mathcal V,\mathbf e}\) be the
canonical diagonals of two tracial matrix ultraproducts.
\begin{enumerate}
\item Every isomorphism
\(\theta:\M_{\mathcal U,\mathbf d}\to\M_{\mathcal V,\mathbf e}\)
sends \(D\) to a unitary conjugate of \(F\).
\item If \(\M_{\mathcal U,\mathbf d}\cong\M_{\mathcal V,\mathbf e}\),
then \((\M_{\mathcal U,\mathbf d},D)\cong(\M_{\mathcal V,\mathbf e},F)\).
\item If \(\M_{\mathcal U,\mathbf d}\cong\M_{\mathcal V,\mathbf e}\),
then
\[
 \prod_\U^{\mathrm{met}}\Sym(d_n)
 \cong\prod_\V^{\mathrm{met}}\Sym(e_n)
\]
as groups with the normalized Hamming metric.
\end{enumerate}
\end{isomorphismtheorem}

\begin{proof}
By part (1), already proved above, we can compose any isomorphism
with an inner automorphism so that it maps \(D\) onto \(F\).
This proves (2). The resulting isomorphism also identifies the
Weyl groups and preserves \(d\).
Using \eqref{eq:weyl} to identify these groups and
\eqref{eq:quotmetric} to identify their metrics, we obtain (3).
\end{proof}

\begin{nonisomorphismtheorem}\label{thm:nonisomorphic-ultraproducts}
There are \(2^{\aleph_0}\) pairwise nonisomorphic tracial matrix
ultraproducts of the sequence \((M_n(\C),\tau_n)\), and
\(2^{\aleph_0}\) pairwise nonisomorphic metric ultraproducts of
\((U(n),d_2)\), where \(d_2(u,v)=\|u-v\|_2\).
\end{nonisomorphismtheorem}

\begin{proof}
Apply \cite[Theorem~5.6 and Remark~5.7]{AT} to choose continuum
many nonisomorphic universal sofic groups. By
Theorem~\ref{thm:isomorphism-reduction}(3), the corresponding
matrix ultraproducts must also be pairwise nonisomorphic.

To obtain the assertion for unitary groups, suppose that two such
groups were isomorphic. Dye's reconstruction theorem would then
give a linear or conjugate-linear \(*\)-isomorphism of the factors;
see \cite[p.~104]{GoldbringSinclair}. Since entrywise complex
conjugation is a conjugate-linear \(*\)-automorphism of a tracial
matrix ultraproduct, either case gives a complex \(*\)-isomorphism,
contradicting the first assertion. Finally, polar decomposition
identifies \(U(\prod_\U M_n)\) with
\(\prod_\U^{\mathrm{met}}(U(n),d_2)\), so these unitary groups
are the required metric ultraproducts.
\end{proof}

In Section~\ref{sec:logical-details} we will express these constructions
by fixed formulas. The diagonal and its normalizer will then be
uniformly definable, and the Weyl group uniformly interpretable.
This will let us apply the same reduction to elementary equivalence.

\section{Internal algebras and isomorphisms}\label{sec:internal-isomorphisms}

To extend diagonal preservation to every internal subalgebra, we
first treat the abelian case. We arrange its coordinate atoms in
rows and columns with distinct total ranks, so that
Lemma~\ref{lem:expander-gaps}(2) gives two pairs of unitaries whose
commutants recover the algebra. An arbitrary internal algebra can
then be recovered inside the commutant of its center using two
further unitaries.

\subsection{Row and column encoding}\label{sec:abelian-encoding}

\begin{proposition}\label{prop:internal-abelian-encoding}
Every internal abelian subalgebra \(B\) of a tracial matrix ultraproduct
\(\M_{\mathcal U,\mathbf d}\) contains commuting internal abelian
algebras \(B_1,B_2\) and admits
four unitaries \(u_1,u_2,v_1,v_2\) such that
\begin{equation}\label{eq:abelian-encoding}
\begin{aligned}
 B_1&=\{u_1,u_2\}',\qquad B_2=\{v_1,v_2\}',\\
 B'&=B_1'\cap B_2',\qquad B=(B_1'\cap B_2')'.
\end{aligned}
\end{equation}
Both pairs \((u_1,u_2)\) and \((v_1,v_2)\) have spectral gap in
\(\M_{\mathcal U,\mathbf d}\) with constant \(8\):
\[
\begin{aligned}
 \|x-E_{B_1}x\|_2^2&\leq8\sum_{j=1}^2\|[x,u_j]\|_2^2,\\
 \|x-E_{B_2}x\|_2^2&\leq8\sum_{j=1}^2\|[x,v_j]\|_2^2
\end{aligned}
\qquad(x\in\M_{\mathcal U,\mathbf d}).
\]
If \(B\) is a MASA, then \(B=B_1'\cap B_2'\).
\end{proposition}

\begin{proof}
We will group the coordinate atoms into rows and columns so that the
block spectral gap estimates apply to both families. By Lemma~\ref{lem:commutant},
we can choose abelian models \(B_n\subseteq M_{d_n}\). We list their atoms \(e_{i,n}\) in
decreasing order of rank and put \(m_{i,n}=\operatorname{rank}(e_{i,n})\), so that
\[
 e_{1,n},\ldots,e_{r_n,n},\qquad
 m_{1,n}\geq\cdots\geq m_{r_n,n}>0.
\]
Our first task is to find projections \(p_n\in B_n'\) with \(\tau_{d_n}(1-p_n)\to_\U0\) and two
abelian algebras \(C_{1,n},C_{2,n}\subseteq p_nB_np_n\), each with pairwise distinct atom ranks,
that together generate \(p_nB_np_n\). The small discarded complements will
not affect the resulting internal algebras.

Suppose first that \(r_n\to_\U\infty\). To arrange the atoms in a triangular
array, let \(t_n=k_n(k_n+1)/2\leq r_n\) be the largest triangular number below \(r_n\)
and retain the first \(t_n\) atoms. We set \(p_n=\sum_{i\leq t_n}e_{i,n}\). Since the
discarded atoms have the smallest ranks, their total trace satisfies
\[
 \tau_{d_n}(1-p_n)
 \leq\frac{r_n-t_n}{r_n}
 \leq\frac{k_n}{r_n}
 \leq\sqrt{\frac2{r_n}}\longrightarrow_\U0.
\]
Thus this choice discards a negligible complement. We arrange the
retained atoms, in their given order, in an array with row lengths
\(k_n,k_n-1,\ldots,1\), filling one row at a time. Suppressing \(n\), we denote
them by \(e_{ij}\), where \(1\leq i\leq k\) and \(1\leq j\leq k+1-i\), and put \(m_{ij}=\operatorname{rank}(e_{ij})\). The
ordering ensures that these ranks are nonincreasing along both rows
and columns.

We now group the atoms by rows and by columns, defining
\[
 E_i=\sum_{j=1}^{k+1-i}e_{ij},\quad
 F_j=\sum_{i=1}^{k+1-j}e_{ij},\quad
 C_1=\bigoplus_{i=1}^k\C E_i,\quad
 C_2=\bigoplus_{j=1}^k\C F_j.
\]
To see why their ranks are distinct, write \(R_i=\operatorname{rank}(E_i)\) and \(S_j=\operatorname{rank}(F_j)\).
For \(1\leq i<k\), the difference of adjacent row ranks is
\[
 R_i-R_{i+1}
 =m_{i,k+1-i}+\sum_{j=1}^{k-i}(m_{ij}-m_{i+1,j})>0,
\]
and the analogous comparison of adjacent columns gives \(S_j>S_{j+1}\) for
\(1\leq j<k\). Hence both families have the required distinct ranks.
Moreover, \(E_iF_j=e_{ij}\) at every occupied cell, so the row and column
algebras together recover \(pB_np\).

If \(r_n\) is instead bounded on a set in \(\U\), we can make the
atom ranks distinct by removing only a bounded number of dimensions.
After restricting to that set, say \(r_n\leq R\), we choose \(e_{i,n}'\leq e_{i,n}\) with
\[
 \operatorname{rank}(e_{i,n}')=\max\{m_{i,n}-i,0\},\qquad
 p_n=\sum_i e_{i,n}'.
\]
The positive retained ranks are now strictly decreasing, and at most
\(R(R+1)/2\) dimensions have been discarded. Consequently \(\tau_{d_n}(1-p_n)\to_\U0\), and
in this case we may simply take \(C_{1,n}=C_{2,n}=p_nB_np_n\).

In either case we obtain internal abelian algebras by setting
\[
 B_i=[C_{i,n}\oplus\C(1-p_n)]_\U\quad(i=1,2).
\]
For \(B_1\), the distinct atom ranks allow us to apply
Lemma~\ref{lem:expander-gaps}(2) to \(C_{1,n}\) in \(p_nM_{d_n}p_n\). We extend the resulting
unitaries by the identity on \(1-p_n\), denote the extensions by
\(u_{1,n},u_{2,n}\), and put \(u_j=[u_{j,n}]_\U\) for \(j=1,2\). If \(p_n=0\), we choose arbitrary
coordinate unitaries; these choices are irrelevant in the ultraproduct
because \(\tau_{d_n}(p_n)\to_\U1\).

It remains to pass the corner estimate to the ambient algebra.
Multiplying it by \(\tau_{d_n}(p_n)\) expresses it in the ambient \(2\)-norm.
Moreover, compression changes a bounded \(x_n\) by at most \(2\|x_n\|\sqrt{\tau_{d_n}(1-p_n)}\)
in \(2\)-norm. Since the unitaries reduce \(p_n\) and the compression
error tends to zero, the estimate from Lemma~\ref{lem:expander-gaps}(2) passes to
the ultraproduct. We obtain \(B_1=\{u_1,u_2\}'\), and \((u_1,u_2)\) has spectral gap
with constant \(8\). Applying the same construction to \(C_{2,n}\)
gives \(v_1,v_2\) with commutant \(B_2\) and the same spectral gap
constant.

We finally verify that these two internal algebras recover the original
one. The unit balls of \(H_n=p_nB_np_n\oplus\C(1-p_n)\) and \(B_n\) have Hausdorff
\(2\)-distance at most \(2\sqrt{\tau_{d_n}(1-p_n)}\). Indeed, compression sends a
contraction in \(B_n\) to a nearby contraction in \(H_n\). In the
other direction, we extend the scalar on each retained part of an
atom to the whole atom, choosing zero on discarded atoms. Since the
discarded trace tends to zero, we have \([H_n]_\U=B\) and \(B_1,B_2\subseteq B\).

The two coordinate algebras commute and generate \(H_n\). To pass
this generation statement to their commutants, let \(P_{i,n}=E_{(C_{i,n}\oplus\C(1-p_n))'}\) for
\(i=1,2\), with commutants taken in \(M_{d_n}\). Since the coordinate algebras commute,
these expectations commute and have product \(E_{H_n'}\). Therefore
Lemma~\ref{lem:commutant}, together with the estimate
\[
 \|x_n-P_{1,n}P_{2,n}x_n\|_2
 \leq\|x_n-P_{1,n}x_n\|_2+\|x_n-P_{2,n}x_n\|_2
\]
gives \(B_1'\cap B_2'=B'\). Taking commutants and using \(B''=B\) proves \eqref{eq:abelian-encoding}.
When the original algebra is a MASA, the final assertion follows
from \(B'=B\).
\end{proof}

\subsection{Preservation of internal algebras}\label{sec:arbitrary-invariance}

The encoding reduces preservation in the abelian case to spectral
gap internality and relative joins of commuting internal abelian
algebras.

\begin{proposition}\label{prop:internal-abelian-invariance}
Every isomorphism \(\theta:\M_{\mathcal U,\mathbf d}\to\M_{\mathcal V,\mathbf e}\)
of tracial matrix ultraproducts maps internal abelian algebras to
internal abelian algebras.
\end{proposition}

\begin{proof}
Let \(\theta\) be such an isomorphism and let \(B\) be internal and
abelian. Since the ambient algebras are finite factors, \(\theta\)
preserves the trace and hence the spectral gap estimates. We choose
\(B_1,B_2\) and \(u_1,u_2,v_1,v_2\) as in Proposition~\ref{prop:internal-abelian-encoding} so that we can transfer
their commutant descriptions to the target. The pair \((\theta(u_1),\theta(u_2))\) then
has spectral gap there, with commutant \(\theta(B_1)\). Theorem~\ref{thm:spectral-gap-internality}
therefore shows that \(\theta(B_1)\) is internal, and the same argument
applies to \(\theta(B_2)\). The two images are abelian and still commute,
so Proposition~\ref{prop:abelian-commuting-internal} makes their relative join internal.
By \eqref{eq:abelian-encoding}, this join is exactly \(\theta(B)\), which proves the assertion.
\end{proof}

For a general internal algebra, we apply the abelian result to its
center. The following relative commutant description then recovers
the algebra inside the commutant of that center.

\begin{lemma}\label{lem:relative-encoding}
Let \(A\) be an internal subalgebra of a tracial matrix ultraproduct
\(\M_{\mathcal U,\mathbf d}\), and put \(Z=Z(A)\).
Then there are \(u_1,u_2\in U(A')\subseteq U(Z')\) such that
\[
 A=Z'\cap\{u_1,u_2\}',
\]
and \((u_1,u_2)\) has spectral gap in \(Z'\) with constant \(2\).
\end{lemma}

\begin{proof}
We write \(A=[A_n]_\U\) and use the multiplicity spaces in the coordinate
decomposition to choose the required unitaries:
\begin{align*}
 A_n&=\bigoplus_i(M_{a_{i,n}}\otimes1_{b_{i,n}}),\\
 u_{j,n}&=\bigoplus_i(1_{a_{i,n}}\otimes v_{j,b_{i,n}}).
\end{align*}
On \(Z(A_n)'=\bigoplus_iM_{a_{i,n}b_{i,n}}\), we can apply Lemma~\ref{lem:expander-gaps}(1) in each multiplicity factor.
Tensoring its estimate with the identity on the first factors and
summing with the inherited trace weights gives the coordinate
spectral gap estimate. By Lemma~\ref{lem:commutant}, the central commutants
have ultraproduct \(Z'\), while the coordinate expectations onto
\(A_n\) induce \(E_A\). Thus both sides of the estimate pass to the
ultraproduct, giving the asserted spectral gap for \(u_j=[u_{j,n}]_\U\in U(A')\).
\end{proof}

\begin{proof}[Proof of Theorem~\ref{thm:invariance}]
Let \(\theta:\M_{\mathcal U,\mathbf d}\to\M_{\mathcal V,\mathbf e}\) be an isomorphism and let \(A\subseteq\M_{\mathcal U,\mathbf d}\) be internal. We first
apply the abelian case to its center \(Z=Z(A)\), which is internal by
Lemma~\ref{lem:commutant}. Proposition~\ref{prop:internal-abelian-invariance} then shows that \(\theta(Z)\) is internal,
and consequently its commutant \(\theta(Z')\) is internal as well.
We can now transfer the description in Lemma~\ref{lem:relative-encoding} to this
internal algebra. It expresses \(\theta(A)\) as the relative commutant
of two unitaries whose spectral gap estimate is preserved by the
isomorphism. Theorem~\ref{thm:spectral-gap-internality} therefore implies that \(\theta(A)\) is
internal. Applying the same argument to \(\theta^{-1}\) proves the converse.
\end{proof}

The same relative description also extends the result on commuting
internal algebras beyond the abelian case.

\begin{corollary}\label{prop:gen-commuting-internal}
If \(A,B\) are commuting internal subalgebras of a tracial matrix
ultraproduct \(\M_{\mathcal U,\mathbf d}\), then
\(A'\cap B'\) and \(A\vee_{\mathrm{rel}}B\) are internal.
\end{corollary}

\begin{proof}
We first put the two commutant descriptions in a common internal
algebra. Since the centers are commuting internal abelian algebras,
\(P_0=Z(A)'\cap Z(B)'\) is internal by Proposition~\ref{prop:abelian-commuting-internal}. Applying Lemma~\ref{lem:relative-encoding} to
\(A'\) and \(B'\), and using \(Z(A')=Z(A)\) and \(Z(B')=Z(B)\), gives \(u_1,u_2\in U(A)\) and
\(v_1,v_2\in U(B)\) that describe \(A'\) and \(B'\) inside their respective
central commutants, with spectral gap constant \(2\). All four
unitaries belong to \(P_0\), where both spectral gap estimates apply.

To do this, consider the expectations \(E_{A'}\) and \(E_{B'}\). They
commute and preserve \(P_0\), since they project onto the invariant
vectors of the commuting conjugation actions of \(U(A)\) and \(U(B)\).
Their product is therefore the expectation onto \(C=A'\cap B'\). Combining
this fact with the two spectral gap estimates gives, for \(x\in P_0\),
\[
 \|x-E_Cx\|_2^2
 \leq\|x-E_{A'}x\|_2^2+\|x-E_{B'}x\|_2^2
 \leq2\sum_{j=1}^2
   \bigl(\|[x,u_j]\|_2^2+\|[x,v_j]\|_2^2\bigr).
\]
It follows that \(u_1,u_2,v_1,v_2\) have spectral gap in \(P_0\), with commutant
\(C=P_0\cap\{u_1,u_2,v_1,v_2\}'\). We can therefore apply Theorem~\ref{thm:spectral-gap-internality} to conclude that
\(C\) is internal. Lemma~\ref{lem:commutant} then gives internality of its
commutant \(C'=A\vee_{\mathrm{rel}}B\) as well.
\end{proof}

\section{Definability of the diagonal and its normalizer}\label{sec:logical-details}
\label{sec:definability-normalizer}

To pass from isomorphisms to elementary equivalence, we express
the constructions of
Sections~\ref{sec:defining-diagonal} and~\ref{sec:internal-isomorphisms}
by uniform distance formulas. Spectral gap defines the required
commutants, and the normalizer estimate then allows us to interpret
the Weyl group. Internality and conjugacy ensure that the diagonal
recovered in this way can be replaced by the canonical one.

\subsection{Distance formulas and internal algebras}\label{sec:distance-formulas}

Write \(\delta_H(x)=\dist_2(x,H_{\leq1})=\|x-E_Hx\|_2\)
for a unital von Neumann subalgebra \(H\) of a tracial algebra \(M\).
All algebra variables below range over \(M_{\leq1}\).
To quantify over a definable set, we replace the restriction by a
distance penalty. For an \(L\)-Lipschitz function \(f\) and a
nonempty subset \(S\) of this ball, we have
\[
 \inf_{s\in S}f(s)=\inf_{\|y\|\leq1}(f(y)+L\dist_2(y,S)),\qquad
 \sup_{s\in S}f(s)=\sup_{\|y\|\leq1}(f(y)-L\dist_2(y,S)).
\]
Both identities follow by taking approximate nearest points.
Thus a formula for the distance to a set lets us replace
quantification over that set by quantification over the whole ball
\cite[Section~3.2]{FHLRTVW}.

\begin{lemma}\label{lem:distance-formulas}
The following distance identities hold in a tracial von Neumann
algebra \(M\).
\begin{enumerate}
\item Let \(P\subseteq M\) be a unital von Neumann subalgebra, and let
\(\mathbf u=(u_1,\ldots,u_h)\in U(P)^h\) have spectral gap in
\(P\) with constant \(L\). For \(C=P\cap\{\mathbf u\}'\),
\[
 \delta_C(x)=\inf_{\|y\|\leq1}\left(\|x-y\|_2+\delta_P(y)
       +\sqrt L\Bigl(\sum_{j=1}^h\|[u_j,y]\|_2^2\Bigr)^{1/2}\right).
\]
\item If \(H_1,\ldots,H_r\subseteq M\) are commuting unital
von Neumann subalgebras and \(K=\bigcap_iH_i'\), then
\[
 \sigma_i(x)=\sup_{\|y\|\leq1}
       (\|[x,y]\|_2\dminus2\delta_{H_i}(y)),\qquad
 \delta_K(x)=\inf_{\|y\|\leq1}
       \Bigl(\|x-y\|_2+\sum_i\sigma_i(y)\Bigr).
\]
\end{enumerate}
In both cases, finite formulas for the input distances give a
finite formula for the resulting distance over the same parameters
and, in (1), the tuple \(\mathbf u\).
\end{lemma}

\begin{proof}
For (1), we first bound the distance to the relative commutant.
Bimodularity of \(E_P\) allows us to apply spectral gap inside
the subalgebra, and Pythagoras then gives
\(\delta_C(x)^2\leq\delta_P(x)^2+L\sum_j\|[u_j,x]\|_2^2\).
Combining this estimate with the triangle inequality shows that
every term in the displayed infimum is at least
\(\delta_C(x)\). Conversely, the choice \(y=E_Cx\) attains this
value, proving the formula; compare
\cite[Lemma~4.1]{GoldbringHartTowsner}.

For (2), the penalty identity first restricts each supremum to the
given algebra, yielding
\(\sigma_i(x)=\sup_{y\in(H_i)_{\leq1}}\|[x,y]\|_2\).
Averaging unitary conjugates then gives
\(\|x-E_{H_i'}x\|_2\leq\sigma_i(x)\), as in
\cite[Corollary~5.12]{Goldbring}. Since the algebras commute,
the expectations onto their commutants commute and have product
\(E_K\). Telescoping this product therefore gives
\(\|x-E_Kx\|_2\leq\sum_i\sigma_i(x)\).
We now argue as in (1): the triangle inequality gives the lower
bound for the infimum, and conditional expectation onto the joint
commutant gives equality.
\end{proof}

\begin{proposition}\label{thm:internal-finite-parameter-definability}
Every internal abelian subalgebra of a tracial matrix ultraproduct
\(\M_{\mathcal U,\mathbf d}\) is definable over at most four unitary
parameters, and every internal subalgebra is
definable over at most six unitary parameters. Their unit-ball
distance predicates are given by fixed finite formulas.
The formulas are independent of the algebra, dimensions, and ultrafilter.
\end{proposition}

\begin{proof}
For an internal abelian algebra \(B\),
Proposition~\ref{prop:internal-abelian-encoding} writes
\(B=(B_1'\cap B_2')'\), where the commuting abelian algebras
\(B_i\) are commutants of unitary pairs with spectral gap
constant \(8\). Lemma~\ref{lem:distance-formulas}(1) therefore
gives their distance formulas over the four parameters. Applying
(2) first to \(B_1,B_2\) and then to their joint commutant gives
the distance to \(B\) without introducing further parameters.

For arbitrary internal \(A\), we apply the abelian construction
to \(Z(A)\). Part (2) then gives the distance to \(P=Z(A)'\)
over the same four parameters. By Lemma~\ref{lem:relative-encoding},
\(A\) is the relative commutant in \(P\) of two further unitaries
with gap constant \(2\), so part (1) gives its distance formula
over all six parameters. Since each step uses fixed formulas and
constants, the resulting formulas have the asserted uniformity.
\end{proof}

\subsection{The diagonal code}\label{sec:diagonal-code}

To transfer a diagonal under elementary equivalence, we need a
condition on the parameters which guarantees that they define an
internal MASA. We obtain it by expressing the conditions of
Lemma~\ref{lem:diagonal-recognition} by finite formulas. The row
and column construction supplies tuples satisfying these conditions.

\begin{lemma}\label{lem:explicit-diagonal-code}
There are fixed nonnegative finite formulas \(\Delta(\mathbf z)\)
and \(\delta(x;\mathbf z)\), with
\(\mathbf z=(z_{ij})_{1\leq i\leq2,\,1\leq j\leq3}\), such that
in a tracial matrix ultraproduct \(M\), a tuple
\(\mathbf z\in M_{\leq1}^6\) satisfies \(\Delta(\mathbf z)=0\)
if and only if:
\begin{enumerate}
\item the \(z_{ij}\) are unitary and their commutants
\(C_i=\{z_{i1},z_{i2},z_{i3}\}'\cap M\) satisfy
\[
 \dist_2(x,(C_i)_{\leq1})\leq4\sum_{j=1}^3\|[z_{ij},x]\|_2
 \qquad(x\in M_{\leq1},\ i=1,2);
\]
\item \(C_1,C_2\) are abelian and commute, and
\(D=C_1'\cap C_2'\) is abelian.
\end{enumerate}
For every such zero, \(D\) is an internal MASA and
\(\delta(x;\mathbf z)=\dist_2(x,D_{\leq1})\).
In \(M_d(\C)\), the same zero-set characterization and distance
identity hold, with \(D\) unitarily conjugate to \(D_d\).
\end{lemma}

\begin{proof}
We encode the distance estimates and the commutation conditions
separately. For \(i=1,2\), define the following formulas on
contractions:
\begin{align*}
 \eta(\mathbf z)&=\max_{i,j}\|z_{ij}^*z_{ij}-1\|_2,
 &\epsilon_i(x)&=\sum_{j=1}^3\|[z_{ij},x]\|_2,\\
 \delta_i(x)&=\inf_{\|y\|\leq1}(\|x-y\|_2+4\epsilon_i(y)),
 &\sigma_i(x)&=\sup_{\|y\|\leq1}
       (\|[x,y]\|_2\dminus2\delta_i(y)),\\
 \gamma_i(\mathbf z)&=\sup_{\|x\|\leq1}\inf_{\|y\|\leq1}
       \max\{\epsilon_i(y),\|x-y\|_2\dminus4\epsilon_i(x)\},
 &\kappa_{ij}(\mathbf z)&=\sup_{\|x\|\leq1}
       (\sigma_i(x)\dminus2\delta_j(x)).
\end{align*}
Here \(\kappa_{ij}\) is indexed by \(1\leq i,j\leq2\).
To recover the joint commutant and impose its abelianness, put
\begin{align*}
 \delta(x;\mathbf z)&=\inf_{\|y\|\leq1}
       (\|x-y\|_2+\sigma_1(y)+\sigma_2(y)),\\
 \beta(\mathbf z)&=\sup_{\|x\|,\|y\|\leq1}
       (\|[x,y]\|_2\dminus2(\delta(x;\mathbf z)+\delta(y;\mathbf z))),\\
 \Delta(\mathbf z)&=\tfrac12\max\{\eta(\mathbf z),
       \max_i\gamma_i(\mathbf z),\max_{i,j}\kappa_{ij}(\mathbf z),\beta(\mathbf z)\}.
\end{align*}

We verify the zero-set characterization in the order of the
construction. The condition \(\eta=0\) gives unitarity, while
compactness in a matrix algebra, or countable saturation in an
ultraproduct, makes \(\gamma_i=0\) equivalent to the distance
estimate in (1). Under this estimate, the triangle inequality
gives \(\delta_i=\delta_{C_i}\); squaring the estimate also gives
spectral gap with constant \(48\). We can therefore use the
penalty identities to identify \(\kappa_{ij}=0\) with commutation
of \(C_i,C_j\). Lemma~\ref{lem:distance-formulas}(2) then gives
the distance formula for the joint commutant, so \(\beta=0\) is
exactly abelianness of \(D\). These identifications prove
both directions of the characterization.

It remains to identify the resulting algebra. By
Lemma~\ref{lem:diagonal-recognition}, \(D\) is an internal MASA
in the ultraproduct case. In a matrix algebra, \(C_i\subseteq D\)
implies \(D'\subseteq D\); since the joint commutant is abelian,
it is a MASA, and simultaneous diagonalization gives the asserted
conjugacy.
\end{proof}

\begin{proof}[Proof of Theorem~\ref{thm:diagonal-code}]
To show that the zero-set is nonempty, use
Lemma~\ref{lem:row-column-recovery} to choose two unitary triples
with spectral gap constant \(10\), whose commutants \(B_1,B_2\)
satisfy \(B_1'\cap B_2'=D\) for the canonical diagonal \(D\).
Their spectral gap estimates imply the required distance bounds,
and hence the triples give a zero of \(\Delta\).
The characterization in Lemma~\ref{lem:explicit-diagonal-code}
applies to every zero, while the conjugacy assertion follows from
Lemma~\ref{lem:commutant}. Together these give the remaining claims.
\end{proof}

\subsection{The Weyl group}\label{sec:normalizer-formulas}

The normalizer estimate turns the distance to the diagonal into
the distance to its normalizer. The same diagonal predicate recovers
the Hamming metric on the quotient, giving a uniform interpretation
of the Weyl group and its action.

\begin{proposition}\label{prop:quotient-interpretation}
Let \(D\) be an internal MASA in a tracial matrix ultraproduct
\(M\), and put \(\delta_D(x)=\dist_2(x,D_{\leq1})\). The formulas
\begin{align*}
 \xi_D(v)&=\sup_{\|y\|\leq1}(\delta_D(vyv^*)\dminus\delta_D(y)),\\
 \nu(x)&=\inf_{\|v\|\leq1}
       \bigl(\|x-v\|_2+\sqrt2\,\xi_D(v)
                +(1+2\sqrt2)\|v^*v-1\|_2\bigr),\\
 d(\bar u,\bar v)&=\delta_D(u^*v)^2\qquad(u,v\in N_M(D))
\end{align*}
give the distance \(\nu(x)=\dist_2(x,N_M(D))\) and the normalized
Hamming metric on \(N_M(D)/U(D)\). Thus this group and its action
on \(D\) are uniformly interpretable in \((M,D)\).
The same assertions hold for \(D_d\subseteq M_d\).
\end{proposition}

\begin{proof}
We first recover the normalizer distance. The penalty identity
restricts the supremum to the diagonal and gives
\(\xi_D(v)=\sup_{y\in D_{\leq1}}\delta_D(vyv^*)\).
For unitary \(v\), this is the normalizer defect in
Lemma~\ref{lem:normalizer}, so its distance to \(N_M(D)\) is at
most \(\sqrt2\,\xi_D(v)\). Since the defect vanishes on the
normalizer, the infimum over unitaries of
\(\|x-v\|_2+\sqrt2\,\xi_D(v)\) is precisely the required distance.
To remove the unitary restriction, note that this integrand is
\((1+2\sqrt2)\)-Lipschitz in \(v\), while polar decomposition gives
\(\dist_2(v,U(M))\leq\|v^*v-1\|_2\) for contractions.
Thus the penalty bounds the change in the integrand when rounding
to a unitary, which proves the formula for \(\nu\).

To identify the quotient metric, consider a coordinate normalizer
\(w=\sigma z\), with \(z\) diagonal unitary. Its diagonal
entries give
\(\delta_{D_d}(w)^2=1-\|E_{D_d}w\|_2^2=d(\sigma,1)\).
By Lemma~\ref{lem:normalizer}, coordinate normalizers represent
every normalizer in the ultraproduct, so this identity passes to
the ultraproduct and gives the formula for \(d\). Distance penalties now interpret the quotient as an
imaginary sort \cite[Section~11]{BBHU}. Its conjugation action
on \(D\) is well defined because \(U(D)\) acts trivially there.
\end{proof}

\subsection{Upgrading to elementary equivalence}\label{sec:elementary-consequences}

We can now transfer the type of a diagonal code between elementarily
equivalent matrix ultraproducts.

\begin{proof}[Proof of Theorem~\ref{thm:reduction}(2)]
Let \(M=\M_{\mathcal U,\mathbf d}\) and
\(P=\M_{\mathcal V,\mathbf e}\) be elementarily equivalent,
with canonical diagonals \(D,F\). We recover the target diagonal
from the type of a tuple coding the source diagonal. By
Theorem~\ref{thm:diagonal-code} and unitary conjugacy, we can
choose \(\mathbf z\) in \(M\) with \(\Delta(\mathbf z)=0\) such
that \(\delta(\,\cdot\,;\mathbf z)\) is the distance to \(D_{\leq1}\).
Elementary equivalence and countable saturation then give a tuple
\(\widetilde{\mathbf z}\) in \(P\) with the same complete type
as \(\mathbf z\).

In particular, \(\Delta(\widetilde{\mathbf z})=0\), so
Theorem~\ref{thm:diagonal-code} identifies
\(\delta(\,\cdot\,;\widetilde{\mathbf z})\) as the distance to an internal
MASA \(\widetilde D\subseteq P\), unitarily conjugate to \(F\).
We can now compare sentences in the expanded pairs: replacing
the diagonal predicate by \(\delta(\,\cdot\,;\mathbf z)\) turns each such sentence
into a formula in the coding tuple. Equality of the tuple types
therefore gives \((M,D)\equiv(P,\widetilde D)\), and conjugating
the target diagonal gives \((M,D)\equiv(P,F)\), as required.
\end{proof}

\begin{proof}[Proof of Theorem~\ref{thm:reduction}(3)]
By part (2), the expanded diagonal pairs are elementarily
equivalent. The uniform interpretation in
Proposition~\ref{prop:quotient-interpretation} therefore gives
elementary equivalence of their Weyl groups with the metric
\(d\). Under the identification \eqref{eq:weyl}, these
are the associated universal sofic groups, proving the assertion.
\end{proof}

\begin{proof}[Proof of Theorem~\ref{cor:many-theories}]
Consider the full sequence \(\Sym(n)\). By
\cite[Theorem~5.6 and Remark~5.7]{AT}, we can choose continuum
many ultrafilters for which the metric ultraproducts have pairwise
distinct theories. By
Theorem~\ref{thm:reduction}(3), elementary equivalence of two
corresponding tracial matrix ultraproducts would imply elementary
equivalence of their universal sofic groups. Hence the matrix
ultraproducts have pairwise distinct theories.

We pass to unitary groups using the reconstruction argument and
the Keisler--Shelah theorem: for tracial matrix ultraproducts,
elementary equivalence of the unitary groups as metric groups
implies elementary equivalence of the algebras; see
\cite[p.~104]{GoldbringSinclair}. Combining this implication
with the polar-decomposition identification used in
Theorem~\ref{thm:nonisomorphic-ultraproducts}, we obtain continuum
many distinct theories of metric ultraproducts of \((U(n),d_2)\).
Since both languages are separable, continuum is also the maximum
possible number of complete theories. Finally, if either sequence
of theories converged, \L{}o\'s' theorem would make all its
nonprincipal ultraproducts elementarily equivalent, contradicting
what we have just proved. Thus neither sequence converges.
\end{proof}

\section{Recovering left-right actions}\label{sec:left-right}

Following the method of \cite[Section~5]{AT}, we realize algebraic
ultraproducts of \(\mathrm{PSL}_2(2^r)\) as definable centralizers
in the Weyl group. Transferring their defining conditions recovers
finite groups whose orders are asymptotic to the matrix dimensions.
Element orders then distinguish an explicit family and rule out
its nontrivial amplifications, proving Theorem~\ref{thm:fundamental-group}.

\subsection{Admissible tuples and their transfer}

Recall that \(G_r=\mathrm{PSL}_2(2^r)\) and
\(|G_r|=2^r(2^{2r}-1)\), for \(r\geq2\).
For a tracial matrix ultraproduct \(M\) with canonical diagonal
\(D\), write \(W=N_M(D)/U(D)\), with normalized Hamming metric
\(d\) and its conjugation action on \(D\). For \(a=(a_1,a_2)\in W^2\),
put
\[
 \mathcal C(a)=\{x\in W:[x,a_1]=[x,a_2]=1\},\qquad
 \epsilon_a(x)=\max_{j=1,2}d(xa_j,a_jx).
\]
Here group commutators mean \([x,y]=xyx^{-1}y^{-1}\).
A group is a \emph{CA-group} if the centralizer of every nonidentity
element is abelian; a subgroup is \emph{\(1\)-discrete} if distinct
elements have distance one.

We impose a distance estimate so that the centralizers admit
restricted quantification, and spectral gap so that they are discrete
and the reconstructed action has a large orbit. Commutation and
the group-theoretic conditions will identify the two regular actions.

\begin{definition}\label{def:admissible}
For constants \(K,L\geq1\), a tuple \((a,t)\in W^2\times W\)
is \emph{admissible} if, putting \(b_j=ta_jt^{-1}\) and
\(c=[b_1,b_2]\), the following hold:
\begin{enumerate}
\item \(\dist(x,\mathcal C(a))\leq L\epsilon_a(x)\) for every \(x\in W\);
\item \(\|f-\tau(f)1\|_2^2\leq K\sum_{j=1}^2\|a_jfa_j^*-f\|_2^2\)
      for every \(f\in D\);
\item \([a_i,b_j]=1\) for \(1\leq i,j\leq2\), and
      \([\mathcal C(a),\mathcal C(b)]=1\);
\item \(\mathcal C(a)\) is a CA-group and
      \(\max_{i=1,2}d([c,b_i c b_i^{-1}],1)=1\).
\end{enumerate}
\end{definition}

The action in (2) is well defined, since any normalizer
representatives of \(a_j\) induce the same action on \(D\).
Conjugation by \(t\) gives (1) and (2) for \(b\), so the two
tuples satisfy the same estimates. Condition (3) also makes
their centralizers mutual:
\(b_j\in\mathcal C(a)\) gives
\(C_W(\mathcal C(a))\subseteq\mathcal C(b)\), and commutation
gives the reverse inclusion.

Expanding generators in the regular actions give such tuples
with constants independent of the dimensions.

\begin{lemma}\label{lem:semantic-models}
There are constants \(K,L\) such that, for every sequence
\(r_j\geq2\) tending to infinity, the canonical diagonal pair
\[
 M=\prod_{\mathcal U}M_{|G_{r_j}|},\qquad
 D=[\ell^\infty(G_{r_j})]_{\mathcal U}
\]
has an admissible tuple with
\[
 \mathcal C(a)=\prod_{\mathcal U}^{\mathrm{alg}}\rho_j(G_{r_j}),
 \qquad
 \mathcal C(b)=\prod_{\mathcal U}^{\mathrm{alg}}\lambda_j(G_{r_j}).
\]
Here \(\lambda(g)x=gx\) and \(\rho(g)x=xg^{-1}\) are the regular actions.
\end{lemma}

\begin{proof}
The construction in \cite[Section~2]{KLN} gives expanding generating
pairs \((s_{1,r},s_{2,r})\) with a uniform spectral gap constant
\(K\). Indeed, its generators are \(c,hch^{-1},h^{-1}ch\),
with \(h^3=1\) and second averaging eigenvalue at most
\(2\sqrt2/3\). The displacement under either conjugate of \(c\)
is at most \(2\|\lambda(h)f-f\|_2+\|\lambda(c)f-f\|_2\),
so \((h,c)\) has uniform spectral gap.
We take \(a_i=[\lambda_j(s_{i,r_j})]_{\mathcal U}\) and let \(t\)
be induced by inversion; then \(b_i=[\rho_j(s_{i,r_j})]_{\mathcal U}\).

To control the distance to the centralizer, take \(u\in\Sym(G_r)\), put
\(a^{(r)}=(\lambda(s_{1,r}),\lambda(s_{2,r}))\) and partition \(G_r\) into
\(E_v=\{g:g^{-1}u(g)=v\}\). Writing \(p_v=|E_v|/|G_r|\),
we apply spectral gap to their indicator functions and obtain
\[
 \dist(u,\rho(G_r))=1-\max_vp_v
 \leq\sum_vp_v(1-p_v)
 \leq2K\sum_{i=1}^2d(u\lambda(s_{i,r}),\lambda(s_{i,r})u)
 \leq4K\epsilon_{a^{(r)}}(u).
\]
Thus we may take \(L=4K\). The estimate shows that permutations
which asymptotically commute with the left generators are
asymptotically equal to right translations. It therefore identifies
the ultraproduct centralizer with the right translations, and
conjugacy gives the left translations. Together with spectral gap,
this proves (1)--(3).

To verify the final condition, recall that the \(G_r\) are nonabelian
simple CA-groups \cite[Proposition~12.2]{ACMM}. For generators \(g_1,g_2\), put
\(c=[g_1,g_2]\ne1\). If both \([c,g_i c g_i^{-1}]=1\), the
CA property makes both generators normalize \(C_{G_r}(c)\),
which would then be a nontrivial proper normal subgroup. Since this
contradicts simplicity, one
commutator is nontrivial, and its regular translation has distance
one from the identity, proving (4).
\end{proof}

Fix these constants \(K,L\) from now on. To transfer the
construction, we express the centralizer distance estimate and
the distance itself by the group formulas
\begin{equation}\label{eq:centralizer-formulas}
 \gamma(a)=\sup_x\inf_y\max\{\epsilon_a(y),d(x,y)\dminus L\epsilon_a(x)\},
 \qquad
 \delta_a(x)=\inf_y(d(x,y)+L\epsilon_a(y))
\end{equation}
whose quantifiers range over \(W\). Vanishing of the first formula
expresses (1); under this condition, the second gives the distance
to the centralizer. Together with discreteness, this permits
transfer of the centralizer and its ordinary group theory along
with the tuple.

\begin{lemma}\label{lem:semantic-transfer}
The centralizers of an admissible tuple are \(1\)-discrete and
uniformly definable, with distance formula \(\delta_a\) for
\(\mathcal C(a)\). If \((M,D)\equiv(P,F)\) are canonical
diagonal pairs of tracial matrix ultraproducts, every admissible
tuple in the first has a counterpart of the same type in the second,
and their centralizers are elementarily equivalent as ordinary groups.
\end{lemma}

\begin{proof}
We prove discreteness by examining fixed points. For
\(u=[u_j]\in\mathcal C(a)\), let \(p_j\) be the diagonal
projection onto the fixed points of a permutation representative
\(u_j\). A fixed point can leave this set only where the permutations
fail to commute, so permutation lifts of the \(a_i\) satisfy
\[
 \|a_{i,j}p_ja_{i,j}^*-p_j\|_2^2
 \leq2d(a_{i,j}u_j,u_ja_{i,j})\longrightarrow_{\mathcal U}0.
\]
Thus \(p=[p_j]\) is fixed by \(a\), and condition (2) gives
\(\tau(p)(1-\tau(p))=0\), so \(d(u,1)=1-\tau(p)\in\{0,1\}\).
Conjugacy proves the same assertion for \(\mathcal C(b)\).

For definability, the distance estimate (1) and the triangle inequality identify
\(\delta_a\) with distance to \(\mathcal C(a)\). Conversely,
in these countably saturated groups, \(\gamma(a)=0\) expresses
(1): it supplies a commuting \(y\) within \(L\epsilon_a(x)\)
of each \(x\). Thus this condition and the definitions of the
centralizers transfer with the tuple type.

By Proposition~\ref{prop:quotient-interpretation}, the Weyl group
and its action on the diagonal are uniformly interpreted in the pair.
Since canonical diagonal pairs are coordinate ultraproducts, they are
countably saturated. We may therefore transfer the complete type
of normalizer representatives of \(a_1,a_2,t\) to \((P,F)\).
Condition (2) transfers as a universal condition on \(D_{\leq1}\),
and the fixed-point argument again gives \(1\)-discreteness.
Distance penalties also permit quantification over the centralizers.
Since on a \(1\)-discrete group ordinary group formulas are continuous
formulas using distance for equality, their theories, including
the CA condition, transfer as well. The same holds for the remaining
commutation and nonmetabelian conditions in (3)--(4).
\end{proof}

\subsection{Finite centralizers and reconstruction}

The cluster argument of \cite[Section~3]{AT} turns a definable,
discrete centralizer into an algebraic ultraproduct of finite groups.
It uses only discreteness and the centralizer distance estimate.
We record this step before using spectral gap on the diagonal to
recover the group orders.

\begin{lemma}\label{lem:finite-centralizer-clusters}
Let \(a=(a_1,a_2)\in(\prod_{\mathcal U}^{\mathrm{met}}\Sym(d_j))^2\),
with permutation representatives \(\mathbf a_j=(a_{1,j},a_{2,j})\).
Suppose that \(\dist(x,\mathcal C(a))\leq L\epsilon_a(x)\)
for every \(x\), and that \(\mathcal C(a)\) is \(1\)-discrete. There are finite
groups \(H_j\) and maps \(\alpha_j:H_j\to\Sym(d_j)\), with
\(\alpha_j(1)=1\), such that
\begin{align*}
 \sup_{g\in H_j}\epsilon_{\mathbf a_j}(\alpha_j(g))&\longrightarrow_{\mathcal U}0,\\
 \sup_{g,h\in H_j}d(\alpha_j(gh),\alpha_j(g)\alpha_j(h))&\longrightarrow_{\mathcal U}0,\\
 \inf_{\substack{g,h\in H_j\\g\ne h}}d(\alpha_j(g),\alpha_j(h))&\longrightarrow_{\mathcal U}1.
\end{align*}
The last infimum is set to one for the trivial group.
The map \([g_j]\mapsto[\alpha_j(g_j)]\) identifies
\(\prod_{\mathcal U}^{\mathrm{alg}}H_j\) with \(\mathcal C(a)\).
\end{lemma}

\begin{proof}
We express the hypotheses by formulas that can be evaluated in the
coordinates. Besides \(\gamma(a)=0\), discreteness gives
\[
 \zeta(a):=\sup_x\bigl(\min\{d(x,1),1-d(x,1)\}\dminus\delta_a(x)\bigr)=0.
\]
We put \(\eta_j=\gamma(\mathbf a_j)+1/j\) and \(v_j=\zeta(\mathbf a_j)\),
which both tend to zero along \(\mathcal U\). These coordinate defects
control correction and separation: for every permutation \(x\),
the definitions supply \(y\) with
\begin{equation}\label{eq:cluster-correction}
 \epsilon_{\mathbf a_j}(y)\leq\eta_j,\qquad
 d(x,y)\leq L\epsilon_{\mathbf a_j}(x)+\eta_j,
\end{equation}
and, since \(\delta_{\mathbf a_j}(x)\leq L\epsilon_{\mathbf a_j}(x)\),
\begin{equation}\label{eq:cluster-separation}
 \min\{d(x,1),1-d(x,1)\}\leq L\epsilon_{\mathbf a_j}(x)+v_j.
\end{equation}
To form the clusters, set \(r=1/(100L)\) and
\(V_j=\{x:\epsilon_{\mathbf a_j}(x)\leq r\}\).
We work on a set in \(\mathcal U\) where \(v_j<1/100\) and
\(\eta_j<r/10\). Inversion invariance and subadditivity of the
defect then give, for \(x,y\in V_j\),
\[
 d(x,y)\leq3/100\quad\text{or}\quad d(x,y)\geq97/100.
\]
Every product of two elements of \(V_j\) can be corrected into
\(V_j\) within \(2Lr+\eta_j<1/5\), by
\eqref{eq:cluster-correction}. As in
\cite[Lemma~4.2]{KunThom}, closeness within \(1/5\) defines an
equivalence relation on \(V_j\), whose classes form a finite group
\(H_j\) by composing representatives and correcting the result.
Indeed, the distance gap makes this multiplication independent of
the choices. The two corrected triple products lie within \(4/5\)
of each other, so the same gap places them in one class and gives associativity.
Identity and inversion are inherited from the symmetric group.

To obtain the required uniform estimates, we correct each
representative once more and choose \(\alpha_j(g)\)
with defect at most \(\eta_j\), taking \(\alpha_j(1)=1\).
Their products belong to \(V_j\) and represent the corresponding
product classes. We may therefore apply
\eqref{eq:cluster-separation}, using the fact that elements of one
class are close and elements of distinct classes are far apart, to obtain
\[
 d(\alpha_j(gh),\alpha_j(g)\alpha_j(h))\leq3L\eta_j+v_j,
 \qquad
 d(\alpha_j(g),\alpha_j(h))\geq1-2L\eta_j-v_j\quad(g\ne h).
\]
These estimates give the three limits and show that the induced
homomorphism is injective. To prove surjectivity, let
\(x=[x_j]\in\mathcal C(a)\). Then \(x_j\in V_j\)
on a set in \(\mathcal U\); its class \(g_j\) satisfies
\(d(x_j,\alpha_j(g_j))\leq L(\epsilon_{\mathbf a_j}(x_j)+\eta_j)+v_j\to_{\mathcal U}0\).
Thus the chosen representatives also give every element of the
centralizer, proving surjectivity. At the unused coordinates we take trivial groups.
\end{proof}

For an admissible tuple, the group-theoretic conditions identify
the finite groups obtained from the cluster lemma. Spectral gap
gives an orbit containing asymptotically all points. The two
commuting faithful actions are then regular there, which recovers
the matrix dimensions.

\begin{proposition}\label{lem:semantic-rigidity}
Let \(M=\M_{\mathcal U,\mathbf d}\) have canonical diagonal
\(D\). If its Weyl group has an admissible tuple \((a,t)\),
then there are integers \(r_j\geq2\) such that
\[
 \frac{|G_{r_j}|}{d_j}\longrightarrow_{\mathcal U}1,
 \qquad
 \mathcal C(a)\cong\prod_{\mathcal U}^{\mathrm{alg}}G_{r_j}.
\]
\end{proposition}

\begin{proof}
Lemmas~\ref{lem:semantic-transfer} and~\ref{lem:finite-centralizer-clusters}
give finite groups \(H_j\) and uniformly separated approximate
representations \(\alpha_j\) realizing \(\mathcal C(a)\).
By admissibility, this group is CA and nonmetabelian: the
commutator in (4) witnesses noncommutation in its derived subgroup.
Since these are first-order conditions, they hold in \(H_j\) on a
set in \(\mathcal U\). The classification of finite CA-groups
identifies the nonmetabelian ones with \(G_r\), \(r\geq2\)
\cite[Proposition~12.2]{ACMM}; see \cite{Weisner,BSW,Suzuki}.
After changing unused coordinates, we may therefore write
\(H_j=G_{r_j}\).

To recover the two commuting actions, choose permutation lifts
\(a_{i,j},t_j\), and put
\(\beta_j(g)=t_j\alpha_j(g)t_j^{-1}\). Their images represent
\(\mathcal C(b)\) and commute asymptotically uniformly with
those of \(\alpha_j\): otherwise maximizing witnesses would
contradict (3) in the ultraproduct. Consequently
\((g,h)\mapsto\alpha_j(g)\beta_j(h)\) is a uniform approximate
action of \(G_{r_j}\times G_{r_j}\). We may therefore apply uniform
flexible permutation stability \cite[Theorem~1.2]{BC} to obtain genuine actions
\[
 \pi_j:G_{r_j}\times G_{r_j}\to\Sym(m_j),\qquad
 m_j\geq d_j,\qquad m_j/d_j\longrightarrow_{\mathcal U}1,
\]
uniformly close to these maps extended by the identity. By the
uniform separation of the original maps, every nonidentity element
of either factor still moves a proportion of points tending uniformly to one.

To find a large orbit, we use the original generators. Since
\(a_i\in\mathcal C(b)\), choose \(h_{i,j}\in G_{r_j}\)
such that \(\pi_j(1,h_{i,j})\) approximates \(a_{i,j}\), extended
by the identity, and
put \(H'_j=\langle h_{1,j},h_{2,j}\rangle\).
We claim that its action has an orbit \(Y_j\) with \(|Y_j|/m_j\to_{\mathcal U}1\).
Indeed, otherwise unions of orbits give invariant sets with
relative sizes in \([\varepsilon/2,1-\varepsilon/2]\) for some
\(0<\varepsilon\leq1/2\): take an orbit in this interval, or
accumulate smaller ones. Intersecting these sets with the original
\(d_j\) points changes their relative sizes by a quantity tending to
zero. Their diagonal projections are asymptotically fixed by \(a_1,a_2\),
since the corrected generators preserve the sets and approximate
the original generators. The resulting ultraproduct projection has
trace bounded away from zero and one, contradicting (2).

Since the first factor commutes with this action, it preserves the
unique largest orbit \(Y_j\).
Its action there is faithful, since a nonidentity element trivial on
\(Y_j\) would move only \(o_{\mathcal U}(m_j)\) points.
The centralizer of a transitive action on \(Y_j\) has at most
\(|Y_j|\) elements, as a commuting permutation is determined
by the image of one point. Faithfulness therefore gives the first
inequality below, while transitivity gives the second and subgroup
inclusion gives the last:
\[
 |G_{r_j}|\leq|Y_j|\leq|H'_j|\leq|G_{r_j}|.
\]
Equality throughout makes both factor actions regular on \(Y_j\), and gives
\(|G_{r_j}|/d_j\to_{\mathcal U}1\), as required.
\end{proof}

\begin{corollary}\label{cor:semantic-dimension-transfer}
Let \(r_j\geq2\), with \(r_j\to\infty\), and suppose that
\[
 \prod_{\mathcal U}M_{|G_{r_j}|}\equiv\M_{\mathcal V,\mathbf e}.
\]
Then there are \(s_j\geq2\) with \(|G_{s_j}|/e_j\to_{\mathcal V}1\),
and \(\prod_{\mathcal U}^{\mathrm{alg}}G_{r_j}\) and
\(\prod_{\mathcal V}^{\mathrm{alg}}G_{s_j}\) are elementarily
equivalent in the ordinary group language.
\end{corollary}

\begin{proof}
Theorem~\ref{thm:reduction}(2) gives elementary equivalence of the
canonical diagonal pairs. We may therefore transfer the tuple from
Lemma~\ref{lem:semantic-models} using Lemma~\ref{lem:semantic-transfer},
and then apply Proposition~\ref{lem:semantic-rigidity} in the target
to obtain the required finite groups and dimension ratios.
\end{proof}

\subsection{Element orders and amplifications}\label{sec:amplifications}

Corollary~\ref{cor:semantic-dimension-transfer} gives both a
dimension ratio and elementary equivalence of the recovered group
ultraproducts. We now use element orders to extract divisibility
information from their theories.

\begin{lemma}\label{lem:psl2-element-orders}
For \(m\geq3\) and \(r\geq2\), the group \(G_r\) contains an
element of order \(2^m-1\) if and only if \(m\mid r\).
\end{lemma}

\begin{proof}
Since the field has characteristic two, \(G_r=\mathrm{SL}_2(2^r)\),
and an element of odd order is diagonalizable over an algebraic closure,
with eigenvalues \(\lambda,\lambda^{-1}\). Frobenius either fixes
or exchanges them, so the order of the element divides
\(2^r-1\) or \(2^r+1\). To determine when \(2^m-1\) divides
either number, write \(r=qm+j\), where \(0\leq j<m\). Modulo \(2^m-1\),
these two integers are \(2^j-1\) and \(2^j+1\). The first is
divisible by \(2^m-1\) precisely when \(j=0\), while
\[
 0<2^j+1\leq2^{m-1}+1<2^m-1.
\]
Conversely, if \(m\mid r\), the cyclic group
\(\mathbb F_{2^r}^{\times}\) contains an element \(\lambda\) of
order \(2^m-1\), and
\(\operatorname{diag}(\lambda,\lambda^{-1})\in G_r\) has that order.
\end{proof}

Existence of an element of any fixed finite order is a first-order
condition in the discrete group language. Thus the lemma and
\L{}o\'s' theorem read divisibility information from the theory of
an algebraic ultraproduct of the \(G_r\): they determine whether
the indices \(r\) are divisible by each fixed \(m\geq3\) along
the ultrafilter.

To apply this information to amplifications, recall that \(Q^t\)
is the normalized corner \(pM_k(Q)p\)
with \((\tau_k\otimes\tau)(p)=t/k\), and that
\(\mathcal F_{\mathrm{fo}}(Q)=\{t>0:Q^t\equiv Q\}\).
For matrix ultraproducts, this compression can be carried out
coordinatewise.

\begin{lemma}\label{lem:coordinate-amplification}
Let \(Q=\M_{\mathcal U,\mathbf d}\) be a tracial matrix
ultraproduct, and let \(t>0\). Then
\[
 Q^t\cong\prod_{\mathcal U}M_{e_i}(\C),\qquad
 e_i=\max\{1,\lfloor t d_i\rfloor\}.
\]
\end{lemma}

\begin{proof}
To realize the amplification coordinatewise, choose an integer
\(k>t\) and projections
\(p_i\in M_{kd_i}(\C)\) of rank \(e_i\). Under the entrywise
identification \(M_k(Q)\cong\prod_{\mathcal U}M_{kd_i}(\C)\),
the projection \(p=[p_i]_{\mathcal U}\) has normalized trace \(t/k\),
so its normalized corner represents the amplification. Coordinate
compression identifies this corner with
\(\prod_{\mathcal U}p_iM_{kd_i}(\C)p_i\cong
\prod_{\mathcal U}M_{e_i}(\C)\). Since
\(e_i/(kd_i)\to_{\mathcal U}t/k>0\), normalizing the corner
traces does not change the ideal of sequences with vanishing
\(2\)-norm.
\end{proof}

The dimension ratio permits only a bounded change in the group
indices. An unbounded set of fixed divisors along the ultrafilter
forces that change to vanish.

\begin{proposition}\label{prop:rigid-dimension-sequences}
Let \(r_i\geq2\), with \(r_i\to\infty\). Suppose
that the integers \(m\geq3\) satisfying
\(\{i:m\mid r_i\}\in\mathcal U\) form an unbounded set. Then
\[
 \mathcal F_{\mathrm{fo}}\left(\prod_{\mathcal U}M_{|G_{r_i}|}(\C)\right)
 =\mathcal F\left(\prod_{\mathcal U}M_{|G_{r_i}|}(\C)\right)=\{1\}.
\]
\end{proposition}

\begin{proof}
Put \(Q=\prod_{\mathcal U}M_{|G_{r_i}|}(\C)\) and suppose that
\(Q^t\equiv Q\). By Lemma~\ref{lem:coordinate-amplification},
\(Q^t\cong\prod_{\mathcal U}M_{e_i}(\C)\), where
\(e_i=\max\{1,\lfloor t |G_{r_i}|\rfloor\}\). We can therefore
apply Corollary~\ref{cor:semantic-dimension-transfer} to obtain integers
\(s_i\geq2\) such that
\[
 \frac{|G_{s_i}|}{e_i}\longrightarrow_{\mathcal U}1,
 \qquad
 \prod_{\mathcal U}^{\mathrm{alg}}G_{s_i}
 \equiv\prod_{\mathcal U}^{\mathrm{alg}}G_{r_i}.
\]
The first limit implies that \(s_i\to_{\mathcal U}\infty\), and
comparison with the original dimensions gives
\[
 \frac{|G_{s_i}|}{|G_{r_i}|}
 =8^{s_i-r_i}\frac{1-4^{-s_i}}{1-4^{-r_i}}
 \longrightarrow_{\mathcal U}t.
\]
Hence \(8^{s_i-r_i}\to_{\mathcal U}t\). Since \(t>0\), the
integers \(s_i-r_i\) are bounded on a set in \(\mathcal U\),
so the ultrafilter selects one of their finitely many values, say \(k\).
Thus \(s_i-r_i=k\) on a set in \(\mathcal U\) and \(t=8^k\).

To show that \(k=0\), fix \(m\geq3\) such that
\(\{i:m\mid r_i\}\in\mathcal U\).
By Lemma~\ref{lem:psl2-element-orders} and \L{}o\'s' theorem,
\(\prod_{\mathcal U}^{\mathrm{alg}}G_{r_i}\) contains an element
of order \(2^m-1\). Elementary equivalence gives such an element
in \(\prod_{\mathcal U}^{\mathrm{alg}}G_{s_i}\), so
\(\{i:m\mid s_i\}\in\mathcal U\). Intersecting these two sets
with \(\{i:s_i-r_i=k\}\) shows that \(m\mid k\). Such \(m\)
are unbounded, so \(k=0\), and consequently \(t=1\). The claim
for the ordinary fundamental group follows by inclusion.
\end{proof}

The indices defining the family in the theorem are eventually
divisible by every fixed power of two. Their odd prime divisors
distinguish the members of the family.

\begin{proof}[Proof of Theorem~\ref{thm:fundamental-group}]
Recall that, for a set \(S\) of odd primes,
\[
 r_i(S)=2^i\prod_{\substack{p\in S\\p\leq i}}p,
 \qquad Q_S=\prod_{\mathcal U}M_{|G_{r_i(S)}|}(\C).
\]
For every \(j\geq2\), the integer \(2^j\) divides \(r_i(S)\)
whenever \(i\geq j\). Proposition~\ref{prop:rigid-dimension-sequences}
therefore gives
\(\mathcal F_{\mathrm{fo}}(Q_S)=\mathcal F(Q_S)=\{1\}\).

To distinguish the theories, suppose that \(Q_S\equiv Q_T\). By
Corollary~\ref{cor:semantic-dimension-transfer}, there are integers
\(s_i\geq2\) such that
\[
 \frac{|G_{s_i}|}{|G_{r_i(T)}|}\longrightarrow_{\mathcal U}1,
 \qquad
 \prod_{\mathcal U}^{\mathrm{alg}}G_{s_i}
 \equiv\prod_{\mathcal U}^{\mathrm{alg}}G_{r_i(S)}.
\]
Since \(|G_{r+1}|/|G_r|>8\), the first limit forces
\(s_i=r_i(T)\) on a set in \(\mathcal U\). For an odd prime \(p\),
Lemma~\ref{lem:psl2-element-orders} now gives
\[
 p\in S
 \ \Longleftrightarrow\
 \{i:p\mid r_i(S)\}\in\mathcal U
 \ \Longleftrightarrow\
 \{i:p\mid r_i(T)\}\in\mathcal U
 \ \Longleftrightarrow\ p\in T.
\]
Thus \(S=T\). Since there are \(2^{\aleph_0}\) subsets of the odd
primes, we obtain the asserted number of theories.
\end{proof}

\subsection*{Acknowledgments}
We thank Ben Hayes, David Jekel, and Sorin Popa for
discussions about tracial matrix ultraproducts, model theory, and
random matrices. A comment by Jekel helped motivate the use of
quantum expanders in the diagonal-recovery argument. We also are indebted to Ilijas Farah whose contributions helped to streamline some of the arguments.

\subsection*{AI assistance}
Large language models were used to develop the proofs and write the
manuscript. They suggested recovering the diagonal by combining
commutants with spectral gap and helped with the subsequent
arguments. The authors directed this work.


\begin{thebibliography}{FHLRA}

\bibitem[AD22]{AD}
D.~Akhtiamov and A.~Dogon,
\emph{On uniform Hilbert Schmidt stability of groups},
Proc.\ Amer.\ Math.\ Soc.\ \textbf{150} (2022), no.~4, 1799--1809.

\bibitem[ALT26]{ALT}
V.~Alekseev, J.~Liu and A.~Thom,
\emph{Spectral gap, internality of relative commutants, and
non-hyperlinear groups},
arXiv:2609.37761 (2026).

\bibitem[AT24]{AT}
V.~Alekseev and A.~Thom,
\emph{On non-isomorphic universal sofic groups},
arXiv:2406.06741 (2024), version~3 (2026).

\bibitem[AG26]{AndoGoldbring}
H.~Ando and I.~Goldbring,
\emph{Matrix ultraproducts have trivial fundamental group},
arXiv:2611.18328 (2026).

\bibitem[ACMM25]{ACMM}
V.~Arvind, P.~J.~Cameron, X.~Ma and N.~V.~Maslova,
\emph{Aspects of the commuting graph}, arXiv:2305.07301, version~5 (2025).

\bibitem[BC23]{BC}
O.~Becker and M.~Chapman,
\emph{Stability of approximate group actions: uniform and probabilistic},
J.\ Eur.\ Math.\ Soc.\ \textbf{25} (2023), no.~9, 3599--3632.

\bibitem[BBHU08]{BBHU}
I.~Ben Yaacov, A.~Berenstein, C.~W.~Henson and A.~Usvyatsov,
\emph{Model theory for metric structures},
in \emph{Model Theory with Applications to Algebra and Analysis},
Vol.~II, London Math.\ Soc.\ Lecture Note Ser.~350,
Cambridge University Press, Cambridge, 2008, 315--427.

\bibitem[BSW58]{BSW}
R.~Brauer, M.~Suzuki and G.~E.~Wall,
\emph{A characterization of the one-dimensional unimodular projective
groups over finite fields}, Illinois J.\ Math.\ \textbf{2} (1958), 718--745.

\bibitem[Far25]{FarahGames}
I.~Farah,
\emph{Games in the matrix: Director's cut},
arXiv:2511.07456 (2025).

\bibitem[FHL\!$^{+}$21]{FHLRTVW}
I.~Farah, B.~Hart, M.~Lupini, L.~Robert, A.~Tikuisis,
A.~Vignati and W.~Winter,
\emph{Model theory of \(C^*\)-algebras},
Mem.\ Amer.\ Math.\ Soc.\ \textbf{271} (2021), no.~1324, viii+127~pp.

\bibitem[FHS13]{FHS1}
I.~Farah, B.~Hart and D.~Sherman,
\emph{Model theory of operator algebras I: stability},
Bull.\ Lond.\ Math.\ Soc.\ \textbf{45} (2013), no.~4, 825--838.

\bibitem[FHS14a]{FHS2}
I.~Farah, B.~Hart and D.~Sherman,
\emph{Model theory of operator algebras II: model theory},
Israel J.\ Math.\ \textbf{201} (2014), no.~1, 477--505.

\bibitem[FHS14b]{FHS3}
I.~Farah, B.~Hart and D.~Sherman,
\emph{Model theory of operator algebras III: elementary equivalence
and \(\mathrm{II}_1\) factors},
Bull.\ Lond.\ Math.\ Soc.\ \textbf{46} (2014), no.~3, 609--628.

\bibitem[FS10]{FarahShelah}
I.~Farah and S.~Shelah,
\emph{A dichotomy for the number of ultrapowers},
J.\ Math.\ Logic \textbf{10} (2010), nos.~1--2, 45--81.

\bibitem[GH01]{GeHadwin}
L.~Ge and D.~Hadwin,
\emph{Ultraproducts of \(C^*\)-algebras},
in \emph{Recent Advances in Operator Theory and Related Topics},
Operator Theory: Advances and Applications~127,
Birkh\"auser, Basel, 2001, 305--326.

\bibitem[Gol23]{Goldbring}
I.~Goldbring,
\emph{Spectral gap and definability},
in J.~Iovino (ed.), \emph{Beyond First Order Model Theory, Volume~II},
Chapman \& Hall/CRC, Boca Raton, FL, 2023, 103--138.

\bibitem[GH23]{GoldbringHartSurvey}
I.~Goldbring and B.~Hart,
\emph{A survey on the model theory of tracial von Neumann algebras},
in I.~Goldbring (ed.), \emph{Model Theory of Operator Algebras},
De Gruyter Series in Logic and Its Applications~11,
De Gruyter, Berlin, 2023, 133--157;
\href{https://arxiv.org/abs/2210.15508}{arXiv:2210.15508}.

\bibitem[GHT18]{GoldbringHartTowsner}
I.~Goldbring, B.~Hart and H.~Towsner,
\emph{Explicit sentences distinguishing McDuff's \(\mathrm{II}_1\) factors},
Israel J.\ Math.\ \textbf{227} (2018), no.~1, 365--377.

\bibitem[GS15]{GoldbringSinclair}
I.~Goldbring and T.~Sinclair,
\emph{Games and elementary equivalence of \(\mathrm{II}_1\)-factors},
Pacific J.\ Math.\ \textbf{278} (2015), no.~1, 103--118.

\bibitem[IJJS26]{IJJS}
V.~Iyer, S.~Jain, S.~Jordan and R.~D.~Somma,
\emph{Efficient quantum circuits for high-dimensional representations
of \(\mathrm{SU}(n)\) and Ramanujan quantum expanders},
preprint (2026), \href{https://arxiv.org/abs/2602.15180v1}{arXiv:2602.15180v1}.

\bibitem[Jek23]{JekelSurvey}
D.~Jekel,
\emph{Free probability and model theory of tracial \(W^*\)-algebras},
in I.~Goldbring (ed.), \emph{Model Theory of Operator Algebras},
De Gruyter Series in Logic and Its Applications~11,
De Gruyter, Berlin, 2023, 215--268.

\bibitem[Jek25]{JekelInformation}
D.~Jekel,
\emph{Information geometry for types in the large-\(n\) limit of
random matrices},
Comm.\ Math.\ Phys.\ \textbf{406} (2025), no.~11, article~272,
67~pp.

\bibitem[Jek26]{JekelProcesses}
D.~Jekel,
\emph{Free information geometry and the model theory of
noncommutative stochastic processes},
arXiv:2604.12212 (2026).

\bibitem[KLN06]{KLN}
M.~Kassabov, A.~Lubotzky and N.~Nikolov,
\emph{Finite simple groups as expanders},
Proc.\ Natl.\ Acad.\ Sci.\ USA \textbf{103} (2006), no.~16, 6116--6119.

\bibitem[KE25]{Kunnawalkam}
S.~Kunnawalkam Elayavalli,
\emph{50 Open Problems: Ultraproduct \(\mathrm{II}_1\) factors},
arXiv:2511.20377 (2025).

\bibitem[KT19]{KunThom}
G.~Kun and A.~Thom,
\emph{Inapproximability of actions and Kazhdan's property (T)},
arXiv:1901.03963 (2019), version~3 (2026).

\bibitem[Li05]{Li}
W.~Li,
\emph{On ultraproducts of operator algebras},
Sci.\ China Ser.~A \textbf{48} (2005), no.~9, 1284--1295.

\bibitem[Pop83]{PopaOrthogonal}
S.~Popa,
\emph{Orthogonal pairs of \(*\)-subalgebras in finite von Neumann
algebras},
J.\ Operator Theory \textbf{9} (1983), no.~2, 253--268.

\bibitem[Pop14a]{PopaKS}
S.~Popa,
\emph{A \(\mathrm{II}_1\) factor approach to the Kadison--Singer
problem},
Comm.\ Math.\ Phys.\ \textbf{332} (2014), no.~1, 379--414.

\bibitem[Pop14b]{PopaIndependence}
S.~Popa,
\emph{Independence properties in subalgebras of ultraproduct
\(\mathrm{II}_1\) factors},
J.\ Funct.\ Anal.\ \textbf{266} (2014), no.~9, 5818--5846.

\bibitem[PSS04]{PopaSinclairSmith}
S.~Popa, A.~M.~Sinclair and R.~R.~Smith,
\emph{Perturbations of subalgebras of type \(\mathrm{II}_1\) factors},
J.\ Funct.\ Anal.\ \textbf{213} (2004), no.~2, 346--379;
corrigendum, J.\ Funct.\ Anal.\ \textbf{235} (2006), 355--356.

\bibitem[Suz57]{Suzuki}
M.~Suzuki, \emph{The nonexistence of a certain type of simple groups
of odd order}, Proc.\ Amer.\ Math.\ Soc.\ \textbf{8} (1957), 686--695.

\bibitem[Tho10]{Thomas}
S.~Thomas,
\emph{On the number of universal sofic groups},
Proc.\ Amer.\ Math.\ Soc.\ \textbf{138} (2010), no.~7, 2585--2590.

\bibitem[Wei25]{Weisner}
L.~Weisner, \emph{Groups in which the normaliser of every element
except identity is abelian}, Bull.\ Amer.\ Math.\ Soc.\ \textbf{31}
(1925), 413--416.

\enlargethispage{\baselineskip}
\end{thebibliography}
\end{document}